\RequirePackage{iftex}
\ifXeTeX
  \input{siam-xetex-compat.tex}
\fi
\documentclass[hidelinks,nohypdvips,onefignum,onetabnum]{siamart251216}
\usepackage{siamStyle}
\ifdefined\pdfoutput\pdfoutput=1\fi

\usepackage[mathscr]{eucal}

\usepackage{bm,nicefrac}
\usepackage[ntheorem]{empheq}
\usepackage{dsfont}
\usepackage{csquotes}
\usepackage{graphicx}
\usepackage{color}

\definecolor{color0}{rgb}{0.12156862745098,0.466666666666667,0.705882352941177}
\definecolor{color1}{rgb}{1,0.498039215686275,0.0549019607843137}
\definecolor{color2}{rgb}{0.172549019607843,0.627450980392157,0.172549019607843}
\definecolor{color3}{rgb}{0.83921568627451,0.152941176470588,0.156862745098039}
\definecolor{color4}{rgb}{0.580392156862745,0.403921568627451,0.741176470588235}
\definecolor{color5}{rgb}{0,0,0}

\definecolor{mycolor1}{rgb}{0.00000,0.44700,0.74100}
\definecolor{mycolor2}{rgb}{0.85000,0.32500,0.09800}
\definecolor{mycolor3}{rgb}{0.92900,0.69400,0.12500}
\definecolor{mycolor4}{rgb}{0.46600,0.67400,0.18800}
\definecolor{mycolor5}{rgb}{0.49400,0.18400,0.55600}
\definecolor{color2}{rgb}{0.172549019607843,0.627450980392157,0.172549019607843}

\newcommand{\zeroVec}{\mathbf{0}}

\newcommand{\R}{\mathbb R}

\newcommand{\dt}{\,\mathrm{d}t}
\newcommand{\ddt}{\tfrac{\mathrm{d}}{\mathrm{d}t}}

\newcommand{\finalTime}{{T}}
\newcommand{\bmu}{\boldsymbol{\mu}}

\newcommand{\switch}{\sigma}
\newcommand{\nrModes}{L}

\newcommand{\nrSwitches}{N}
\newcommand{\switchingSet}{\calK}

\newcommand{\system}{\Sigma}
\newcommand{\switchedSys}{\system_{\switch}}
\newcommand{\switchedSysAdj}{\system_{\switch}^\T}
\newcommand{\state}{\mathbf{x}}
\newcommand{\stateDim}{n}
\newcommand{\adState}{\mathbf{p}}
\newcommand{\inp}{\mathbf{u}}
\newcommand{\inpDim}{m}
\newcommand{\out}{\mathbf{y}}
\newcommand{\outDim}{p}
\newcommand{\adOut}{\mathbf{q}}
\newcommand{\mM}{\mathbf{M}}
\newcommand{\mN}{\mathbf{N}}
\newcommand{\mK}{{\mathbf{K}}}
\newcommand{\mE}{\mathbf{E}}
\newcommand{\mA}{\mathbf{A}}
\newcommand{\mC}{\mathbf{C}}
\newcommand{\mB}{\mathbf{B}}

\newcommand{\mS}{\mathbf{S}}
\newcommand{\mR}{\mathbf{R}}
\newcommand{\mW}{\mathbf{W}}
\newcommand{\fU}{\mathbf{U}}

\newcommand{\mI}{\mathbf{I}}
\newcommand{\mV}{\mathbf{V}}
\newcommand{\mx}{\mathbf{x}}

\newcommand{\mv}{\mathbf{v}}
\newcommand{\mPsi}{\mathbf{\Psi}}

\newcommand{\Solper}{\calG}

\newcommand{\reduce}[1]{\tilde{#1}}
\newcommand{\reduceAd}[1]{\hat{#1}}
\newcommand{\systemRed}{\reduce{\system}}
\newcommand{\switchedSysRed}{\systemRed_{\switch}}
\newcommand{\stateDimRed}{r}
\newcommand{\stateRed}{\reduce{\state}}
\newcommand{\outRed}{\reduce{\out}}
\newcommand{\adOutRed}{\reduce{\adOut}}
\newcommand{\adOutRedAd}{\reduceAd{\adOut}}
\newcommand{\adStateRed}{\reduce{\adState}}
\newcommand{\adStateRedAd}{\reduceAd{\adState}}

\newcommand{\mERed}{\reduce{\mathbf{E}}}
\newcommand{\mARed}{\reduce{\mathbf{A}}}
\newcommand{\mCRed}{\reduce{\mathbf{C}}}
\newcommand{\mBRed}{\reduce{\mathbf{B}}}

\newcommand{\cost}{J}

\newcommand{\desired}[1]{#1_{\mathrm{d}}}
\newcommand{\inpVarDes}{\desired{\inp}}
\newcommand{\outVarDes}{\desired{\out}}
\newcommand{\opt}[1]{#1^{\star}}
\newcommand{\inpOpt}{\opt{\inp}}
\newcommand{\outOpt}{\opt{\out}}
\newcommand{\adStateOpt}{\opt{\adState}}
\newcommand{\dadStateOpt}{\opt{\dot{\adState}}}
\newcommand{\stateOpt}{\opt{\state}}
\newcommand{\dstateOpt}{\opt{\dot{\state}}}
\newcommand{\adOutOpt}{\opt{\adOut}}

\newcommand{\multiplier}{\boldsymbol{\xi}}
\newcommand{\multiplierOpt}{\opt{\multiplier}}

\newcommand{\inpRedOpt}{\opt{\reduce{\inp}}}
\newcommand{\outRedOpt}{\opt{\outRed}}
\newcommand{\stateRedOpt}{\opt{\stateRed}}
\newcommand{\dstateRedOpt}{\opt{\dot{\stateRed}}}
\newcommand{\adStateRedOpt}{\opt{\adStateRed}}
\newcommand{\dadStateRedOpt}{\opt{\dot{\adStateRed}}}
\newcommand{\adOutRedOpt}{\opt{\adOutRed}}
\newcommand{\multiplierRedOpt}{\opt{\reduce{\multiplier}}}

\newcommand{\StateVar}{\theta}
\newcommand{\bx}{\mathrm{x}}
\newcommand{\bs}{\mathrm{s}}

\newcommand{\reg}{\gamma}
\newcommand{\Proj}{P_{\mYad}}

\newcommand{\shift}{\mathbf{s}}

\newcommand{\mU}{{\mathscr U}}
\newcommand{\mY}{{\mathscr Y}}

\newcommand{\mYad}{{\mY_\mathsf{ad}}}
\newcommand{\mUad}{{\mU_\mathsf{ad}}}
\newcommand{\Uad}{U_{\mathsf{ad}}}
\newcommand{\Yad}{Y_{\mathsf{ad}}}

\newcommand{\calA}{\mathcal{A}}

\newcommand{\calC}{\mathcal{C}}

\newcommand{\calG}{\mathcal{G}}

\newcommand{\calJ}{\cost}
\newcommand{\calK}{\mathcal{K}}
\newcommand{\calL}{\mathcal{L}}
\newcommand{\calM}{\mathcal{M}}

\newcommand{\calO}{\mathcal{O}}
\newcommand{\calP}{\mathcal{P}}
\newcommand{\calQ}{\mathcal{Q}}
\newcommand{\calR}{\mathcal{R}}

\newcommand{\abbr}[1]{\textsf{\upshape #1}\xspace}

\newcommand{\MY}{\abbr{MY}}

\newcommand{\FOM}{\abbr{FOM}}
\newcommand{\ROM}{\abbr{ROM}}
\newcommand{\ROMs}{\abbr{ROMs}}

\newcommand{\MOR}{\abbr{MOR}}

\newcommand{\OCP}{\abbr{OCP}}
\newcommand{\SLS}{\abbr{SLS}}
\newcommand{\OCPs}{\abbr{OCPs}}
\newcommand{\BT}{\abbr{BT}} 
\newcommand{\RB}{\abbr{RB}} 
\newcommand{\POD}{\abbr{POD}} 

\newcommand{\LTI}{\abbr{LTI}}
\newcommand{\LMI}{\abbr{LMI}}
\newcommand{\LMIs}{\abbr{LMIs}}

\newcommand{\GLEs}{\abbr{GLEs}}
\newcommand{\PBR}{\abbr{PBR}}

\newcommand{\fcP}{\mathbf{\calP}}
\newcommand{\fcQ}{\mathbf{\calQ}}
\newcommand{\T}{\ensuremath\mathsf{T}}
\newcommand{\lyapOper}{\calL}
\newcommand{\gleOper}{\Pi}

\DeclareMathOperator{\vec2}{vec}

\headers{Model reduction for output-constrained control}
{M. Kartmann, M. Manucci, B. Unger, and S. Volkwein}
\title{Balancing-based model reduction for output-constrained optimal control of switched linear systems\thanks{%
\textbf{Funding:} All authors were supported by the BMBF (grant no.~05M22VSA).}}
\author{Michael Kartmann\thanks{FB Mathematik und Statistik, Universität Konstanz,
Universitätsstr. 10, 78464 Konstanz, Germany
(\email{michael.kartmann@uni-konstanz.de}, \email{stefan.volkwein@uni-konstanz.de}).}
\and Mattia Manucci\thanks{Institute of Applied Analysis and Numerical Simulation,
University of Stuttgart, Pfaffenwaldring 57, 70569 Stuttgart, Germany
(\email{mattia.manucci@mathematik.uni-stuttgart.de}).}
\and Benjamin Unger\thanks{Institute for Applied and Numerical Mathematics,
Karlsruhe Institute of Technology, 76131 Karlsruhe, Germany
(\email{benjamin.unger@kit.edu}). This author was additionally supported in part by the Deutsche Forschungsgemeinschaft (DFG, German Research Foundation) -- Project-ID 258734477 -- SFB 1173.}
\and Stefan Volkwein\footnotemark[2]}

\begin{document}

\maketitle

\begin{abstract}
    We develop a balanced-truncation (BT) framework for model-order reduction of optimal control problems governed by switched linear systems, with quadratic tracking costs, convex input constraints, and Moreau--Yosida-regularized output constraints. For general Petrov--Galerkin reduced models, we derive an a posteriori bound for the error in the optimal value function, and combine it with BT-based estimates of the dynamical-system reduction error. We show that BT preserves the duality structure underlying the first-order optimality conditions, making first-reduce-then-optimize and first-optimize-then-reduce equivalent. Consequently, reducing the forward system also reduces the associated optimality system, leading to a certified, online-efficient, and online-adaptive approximation of the parameterized solution operator. The reduced model can therefore be reused across different regularization parameters, constraint sets, switching signals, and reference trajectories.
\end{abstract}

\begin{keywords}
Switched linear systems, model reduction, balanced truncation, optimal control.
\end{keywords}

\begin{MSCcodes}
49K20, 49M05, 93A15, 93C30
\end{MSCcodes}


\section{Introduction}\label{sec:intro}
We consider the optimal control of \emph{switched linear systems}~(\SLS) of the form
\begin{equation}
    \label{eq:SLS:FOM}
    \switchedSys \quad \left\{\quad \begin{aligned}
        \mE_{\switch(t)}\dot{\state}(t)  &= \mA_{\switch(t)} \state(t) +\mB_{\switch(t)}\inp(t), & \state(0) &= \zeroVec, \\
        \state(t^{+})& =\mK_{\switch(t^+),\switch(t^-)}\mx(t^-),\\
        \out(t) &= \mC_{\switch(t)}\state(t),\\
    \end{aligned}\right.
\end{equation}
where $\switch\colon \R\to \switchingSet\coloneqq \{1,\ldots,\nrModes\}$ is the external switching signal, which we assume to be an element of the set
\begin{equation}
    \label{eq:suitableSwitchingSignals}
    \calP \;\vcentcolon=\; \{\switch\colon \R\to \switchingSet \mid \switch \text{ is left-continuous with locally finitely many jumps}\}.
\end{equation}
The symbols $\state(t)\in\R^{\stateDim}$, $\inp(t)\in\R^{\inpDim}$, and $\out(t)\in\R^{\outDim}$ denote the \emph{state}, \emph{control} (or \emph{input}), and measured \emph{output}, respectively.
We use $\state(\inp)$ and $\out(\inp)$ to emphasize their dependence on the control $\inp\in L^2(0,T;\R^\inpDim)$ with terminal time $\finalTime> 0$.
For each mode $j\in\switchingSet$, the matrices $\mE_j,\mA_j\in\R^{\stateDim\times\stateDim}$, $\mB_j\in\R^{\stateDim\times\inpDim}$, and $\mC_j\in\R^{\outDim\times\stateDim}$ define the active subsystem and $\mE_j$ is assumed to be invertible.
The matrices $\mK_{i,j}\in\R^{\stateDim\times\stateDim}$ describe state transitions from mode $j$ to mode $i$ at the switching times, with $\mK_{i,i}=\mI_{\stateDim}$.
As defined in \Cref{sec:notation}, $\switch(t^-)$ and $\switch(t^+)$ denote the left and right limits, respectively.
We refer to~\eqref{eq:SLS:FOM} as the \emph{full-order model} (\FOM).
For each fixed $\switch\in\calP$, let $0=t_0<\cdots<t_{\nrSwitches+1}=\finalTime$ be the partition induced by the switching times, i.e., the discontinuities of $\switch$ in the interval $(0,\finalTime)$.
The state and output satisfy $\state|_{(t_k,t_{k+1})}\in H^1(t_k,t_{k+1};\R^\stateDim)\hookrightarrow C([t_k,t_{k+1}];\R^\stateDim)$ and $\out|_{(t_k,t_{k+1})}\in C([t_k,t_{k+1}];\R^\outDim)$ for $k=0,\ldots,\nrSwitches$.

Given a cost parameter $\lambda>0$, we consider the finite-horizon optimal control of~\eqref{eq:SLS:FOM} with respect to the quadratic cost function
\begin{align}\label{eqn:def:cost:func}
    \begin{aligned}
        \cost(\out, \inp)\;&\vcentcolon=\; \int_{0}^{\finalTime}\left( \tfrac{1}{2}\|\out(t)-\outVarDes(t)\|^2_{\R^\outDim}+\tfrac{\lambda}{2}\|\inp(t)-\inpVarDes(t)\|^2_{\R^\inpDim}\right)\dt,
    \end{aligned}
\end{align}
where $\inpVarDes\in L^2(0,\finalTime;\R^\inpDim)$ and $\outVarDes\in L^2(0,\finalTime;\R^\outDim)$ are prescribed references.
Let $\Uad\subseteq\R^\inpDim$ and $\Yad\subseteq\R^\outDim$ be nonempty, closed, convex sets of admissible inputs and outputs, respectively, and define $\mUad\coloneqq L^2(0,\finalTime;\Uad)$.
The pointwise output-constrained \emph{full-order optimal control problem}~(\FOM-\OCP) is
\begin{equation}\label{eq:OPT:FOMoutcons}
    \min \cost(\out,\inp)\;\; \text{s.t.}\;\;(\out, \inp)\text{ solves }\eqref{eq:SLS:FOM},\;\; \inp\in\mUad,\;\; \out(t)\in\Yad\text{ for all }t\in[0,T].
\end{equation}
By piecewise continuity and closedness of $\Yad$, the output constraints in \eqref{eq:OPT:FOMoutcons} are equivalent to $\out(t)\in\Yad$ for almost every $t\in(0,\finalTime)$, in particular, both one-sided limits at each interior switching time belong to $\Yad$.
A natural setting for multiplier theory is the product of continuous-function spaces on the closed switching intervals, in which Lagrange multipliers consist of tuples of finite vector-valued measures (cf., e.g., \cite[Sec.~6.2.1]{troltzsch2024optimal}).
To obtain computable multipliers, one can consider a \emph{Moreau--Yosida}~(\MY) regularization in $L^2$ (cf.~\cite{bauschke2011convex}, \cite[Sec.~2.7.2.1]{HPUU09}).
To this end, define the nonempty, closed, convex set
\begin{equation*}
    \mYad\vcentcolon=\bigl\{\mv\in L^2(0,T;\R^\outDim)\,\big|\,\mv(t)\in\Yad\text{ for almost every }t\in(0,T)\bigr\}.
\end{equation*}
For a regularization parameter $\reg>0$, the \MY approximation is based on the metric projection
\begin{equation*}
    \Proj\colon L^2(0,T;\R^\outDim)\to\mYad,
    \qquad
    \Proj(\out)\vcentcolon=\operatorname*{arg\,min}_{\mv\in\mYad}\|\out-\mv\|_{L^2(0,T;\R^\outDim)}.
\end{equation*}
With the regularized cost functional $\cost_\reg(\out,\inp)\vcentcolon=\cost(\out,\inp)+\tfrac{\reg}{2}\|\out-\Proj (\out)\|_{L^2(0,T;\R^\outDim)}^2$ we obtain the \MY-regularized \FOM-\OCP
\begin{equation}\label{eq:OPT:FOMoutconsMY}
    \min \cost_\reg(\out,\inp) \quad \text{s.t.}\quad(\out, \inp)\text{ solves }\eqref{eq:SLS:FOM} \text{ and }\inp\in\mUad.
\end{equation}
We denote the unique solution of \eqref{eq:OPT:FOMoutconsMY} by $(\outOpt,\inpOpt)$; see \Cref{rem:FOM:opsysoutconsMY}. The convergence of solutions of~\eqref{eq:OPT:FOMoutconsMY} towards the solution of \eqref{eq:OPT:FOMoutcons} as $\reg\to \infty$ is investigated, e.g., in \cite{ITO2003221}. Approximating the solution of \eqref{eq:OPT:FOMoutconsMY} numerically requires repeated evaluations of \eqref{eq:SLS:FOM} for different inputs $\inp\in\mUad$.
These evaluations become computationally expensive when $\stateDim$ is large, particularly if the \OCP must be solved for many configurations of the parameter tuple
\begin{equation}\label{eq:definition_mu}
    \bmu \vcentcolon= (\sigma,\outVarDes,\inpVarDes,\Uad ,\Yad,\lambda, \reg,\finalTime).
\end{equation}

Our objective is, therefore, to develop a rapidly computable approximation of the \OCP\ solution operator
\begin{equation}\label{eq:sol_op}
    \Solper\colon\bmu\mapsto (\outOpt(\bmu),\inpOpt(\bmu)).
\end{equation}
To construct this approximation, we apply projection-based \emph{model-order reduction} (\MOR) and replace~\eqref{eq:SLS:FOM} by the \emph{reduced-order model} (\ROM)
\begin{equation}
    \label{eq:SLS:ROM}
    \switchedSysRed \quad \left\{\quad \begin{aligned}
        \reduce{\mE}_{\switch(t)}\dot{\stateRed}(t)  &= \reduce{\mA}_{\switch(t)} \stateRed(t) + \reduce{\mB}_{\switch(t)}\inp(t), &    \stateRed(0) &= \zeroVec, \\
        \stateRed(t^{+})& =\reduce{\mK}_{\switch(t^+),\switch(t^-)}\stateRed(t^-),\\
        \reduce{\out}(t) &= \reduce{\mC}_{\switch(t)}\stateRed(t),\\
    \end{aligned}\right.
\end{equation}
with the reduced matrices $\reduce{\mE}_{j},\reduce{\mA}_{j},\reduce{\mK}_{i,j}\in\R^{\stateDimRed\times \stateDimRed}$, $\reduce{\mB}_{j}\in \R^{\stateDimRed\times \inpDim}$, $\reduce{\mC}_{j}\in\R^{\outDim\times \stateDimRed}$, and $\stateDimRed\ll \stateDim$.
The reduced system matrices are obtained via Petrov--Galerkin projection.
We choose full-column-rank matrices $\mV_j,\mW_j\in\R^{\stateDim\times\stateDimRed}$ such that $\mW_j^\T\mE_j\mV_j$ is invertible for each mode $j$, and define
\begin{align}
    \label{eqn:SLS:ROM:matrices}
    \begin{gathered}
        \reduce{\mE}_j  \vcentcolon= \mW_j^\T \mE_j \mV_j, \qquad
        \reduce{\mA}_j  \vcentcolon= \mW_j^\T \mA_j \mV_j, \qquad
        \reduce{\mB}_j  \vcentcolon= \mW_j^\T \mB_j, \\
        \reduce{\mC}_j  \vcentcolon= \mC_j \mV_j, \qquad
        \reduce{\mK}_{i,j} \vcentcolon= \reduce{\mE}_i^{-1}\mW_i^\T\mE_i\mK_{i,j}\mV_j.
    \end{gathered}
\end{align}
The transition formula follows by imposing $\mW_i^\T\mE_i(\mV_i\stateRed(t^+)-\mK_{i,j}\mV_j\stateRed(t^-))=\zeroVec$ at a switch from mode $j$ to mode $i$. By construction, we obtain $\reduce{\mK}_{i,i}=\mI_{\stateDimRed}$.
%
The \MY-regularized \ROM-\OCP reads as
\begin{equation}\label{eq:OPT:ROMoutconsMY}
    \min \cost_\reg(\outRed,\inp)\quad \text{s.t.}\quad(\outRed, \inp)\text{ solves }\eqref{eq:SLS:ROM} \text{ and }\inp\in\mUad,
\end{equation}
and, for each configuration $\bmu$, \eqref{eq:OPT:ROMoutconsMY} has the unique solution $(\outRedOpt(\bmu),\inpRedOpt(\bmu))$.
We do not study convergence as $\reg\to\infty$ nor do we quantify the regularization error compared to \eqref{eq:OPT:FOMoutcons}. Instead, for each prescribed $\reg>0$, we certify the \MOR error between the \FOM-\OCP and the \ROM-\OCP for a fixed $\reg$.
Thus, the goal of certified \MOR is to seek projection matrices $\mW_j$ and $\mV_j$ for which the \ROM-\OCP \eqref{eq:OPT:ROMoutconsMY} approximates the \FOM-\OCP \eqref{eq:OPT:FOMoutconsMY} accurately for all configurations $\bmu$, and a rigorous, possibly online-efficient error bound $\Delta_{\outOpt,\inpOpt}(\outRedOpt(\bmu),\inpRedOpt(\bmu))$ satisfying
\begin{equation}\label{eq:BT_optimal_bound}
    \|(\outRedOpt(\bmu),\inpRedOpt(\bmu))-(\outOpt(\bmu),\inpOpt(\bmu))\| \leq \Delta_{\outOpt,\inpOpt}(\outRedOpt(\bmu),\inpRedOpt(\bmu))\quad \text{ for all } \bmu,
\end{equation}
where $\|\cdot\|$ denotes a suitable norm.
By \emph{online-efficient}, we mean that the computational cost required to evaluate the bound is independent of the full-order dimension $\stateDim$.
Expressing the optimal input and output as functions of $\bmu$ highlights the requirement that $\bmu$ can be varied.

Systems of the form \eqref{eq:SLS:FOM} arise in engineering applications that operate in several modes or regimes.
Optimal control must ensure efficient, safe, and reliable operation while respecting input and output constraints.
The combination of switching dynamics, constraints, and real-time requirements makes the resulting problems computationally demanding and motivates our study.
\subsection{Literature review}
\MOR is a highly active research field that has seen tremendous development in recent years.
For an introduction and overview, we refer to the monographs and collections \cite{BenGW15,BOCW17,HRS16,QMN16}.
Regarding \MOR for \SLS, one straightforward strategy is to construct a separate \ROM for each subsystem \cite{PapP14,PapP16}.
Because it ignores transitions between modes, however, this strategy can produce arbitrarily large output errors; see \cite[Sec.~2]{SchU18}.
Coupled \emph{linear matrix inequalities}~(\LMIs) offer an alternative \cite{WuZ09}, but become computationally infeasible for large-scale systems.
Related matrix-inequality techniques yield balancing-based \ROMs with theoretical guarantees of quadratic stability and an associated error bound \cite{PetWL13}.

The approaches in \cite{PonGB20,SchU18} instead reinterpret the switching signal as an input and reformulate the \SLS as a non-switched system.
The method in \cite{SchU18} embeds the dynamics into a larger \emph{linear time-invariant} (\LTI) system, whereas \cite{PonGB20} uses a bilinear reformulation and the corresponding balanced-truncation method.
Both approaches provide error guarantees, but only under restrictive assumptions.
The \emph{piecewise balanced reduction} (\PBR) method for \SLS introduced in \cite{ManU26}, and subsequently extended to handle switched systems of linear differential-algebraic equations \cite{ManU26b}, relaxes these assumptions by constructing projection matrices from approximate solutions of suitably defined \emph{generalized Lyapunov equations} (\GLEs).
It provides an a posteriori error estimator under general assumptions, although the efficiency of its evaluation depends on the given \SLS.
Because \PBR is the specific \MOR method we chose to efficiently approximate the \OCP, we review it in more detail in the forthcoming \Cref{sec:BT_for_OCP}.
Interpolation-based alternatives for hybrid and switched systems are discussed in~\cite{BasPWL16,ScaA16}.

In the context of optimal control, \MOR is used, e.g., in \cite{Ant09,BBH18,GV17,HJ18,HV08,KarMUV26,NMT11,SBMR18,TV09}.
\emph{Balanced truncation}~(\BT) and related balancing-based methods have been used for controller reduction and feedback design; see, e.g., \cite{BMS06,breiten2015feedback,enns1984model,MusG91}.
Certified \emph{reduced basis}~(\RB) approaches to parametrized control-constrained \OCPs are studied in, e.g., \cite{bader2016certified,DieG23}.
The contributions \cite{Bader674368,gubisch2014pod} also discuss regularization-based approaches to state constraints and a posteriori error estimation with \RB and 
\emph{proper orthogonal decomposition} (\POD) models, respectively.
\subsection{Outline and main contributions}
%
We combine output-based a posteriori error analysis with \PBR to approximate the solution operator \eqref{eq:sol_op} for \MY-regularized optimal control of \SLS.
Our results quantify the model-reduction error at a prescribed regularization parameter and consist of three contributions.

First, for piecewise Petrov--Galerkin \ROMs, \Cref{sec:outputbased_error_est} establishes a joint error bound for the optimal control, output, and approximate multiplier (\Cref{lemma:apost_control}) and two-sided bounds for the difference of the optimal value functions (\Cref{lemma:apost_valuefun}).
These estimates use errors in the forward and adjoint outputs, which allows them to be combined with input-output error bounds from balancing-based reduction.
They require no full-order optimization; their direct evaluation uses forward and adjoint \FOM simulations driven by reduced optimal quantities.
We therefore call them \emph{full-order-evaluation bounds}.

Second, building on the \PBR method of \cite{ManU26}, we prove that reduction commutes with the derivation of the optimality system (\Cref{teo:duality}).
For asymptotically stable modes and uniquely solvable \GLEs~(see \Cref{ass:BT,ass:uniq:GLE}), the reachability and observability matrices of the forward system and its time-reversed adjoint are interchanged.
The adjoint of the reduced forward system consequently coincides with the system obtained by reducing the full-order adjoint.
Thus, the forward balancing construction also determines the reduced adjoint, without a separate balancing computation.

Third, \Cref{cor:PBRcontrol} combines certified forward and adjoint \PBR error estimates with the output-based bounds to certify the optimal control, output, approximate multiplier, and value function.
\Cref{subsec:computational_operator_approx} uses these estimates to select the reduced dimension adaptively from a precomputed hierarchy.
Online efficiency requires a certified offline--online decomposition of the estimator, e.g.,~\eqref{eq:est_cheap} gives an explicit realization when the singular-value-tail bound \eqref{eqn:PBR:special:bound} applies.
When this bound holds uniformly over the bounded parameter set specified in \Cref{rem:uniform:solution}, that corollary proves uniform approximation of the solution operator. For fixed system matrices, the projection bases can thus be reused as the regularization parameters, cost weights, constraint sets, switching signal, and reference trajectories vary.

Finally, \Cref{sec:Numerics} studies a synthetic switched system and a spatially discretized two-room heating problem.
The numerical experiments display approximation errors, estimator sharpness, online computational costs, and computed control and output trajectories.
%
\subsection{Notation}\label{sec:notation}
%
We define the left and right limits as $\switch(t^-)\coloneqq \lim_{s\uparrow t}\switch(s)$ and $\switch(t^+)\coloneqq \lim_{s\downarrow t}\switch(s)$, if they exist.
If it is clear from the context, we abbreviate, e.g., $L^2(0,\finalTime;\R^\inpDim)$ or $L^2(0,\finalTime;\R^\outDim)$ as $L^2$.
For $\state, \adState\in L^2(0,\finalTime;\R^\stateDim)$, we abbreviate, e.g., $\langle  \mA_{\switch}\state, \adState\rangle_{\R^\stateDim} (t)\coloneqq \langle \mA_{\switch(t)}\state(t), \adState(t) \rangle_{\R^\stateDim} $ for almost all $t\in (0,\finalTime)$.
\section{Error estimates for the optimal control, output, and value function}\label{sec:outputbased_error_est}
%
Throughout this section, we fix a configuration $\bmu$ as in \eqref{eq:definition_mu} and compare the optimal solutions of the \OCPs \eqref{eq:OPT:FOMoutconsMY} and \eqref{eq:OPT:ROMoutconsMY}.
We first introduce the adjoint system and the regularized multiplier, then state the optimality conditions on which the error analysis is based.
The gradient of the \MY penalty term in \eqref{eq:OPT:FOMoutconsMY} is the approximate multiplier $\multiplier(\out)=\reg(\out-\Proj(\out))\in L^2(0,\finalTime;\R^\outDim)$.
For arbitrary $\out,\multiplier\in L^2(0,\finalTime;\R^\outDim)$, we denote by $\adState(\out,\multiplier)$ and $\adOut(\out,\multiplier)$ the state and output of
\begin{equation}
    \label{eq:SLS:FOM_Adj}
    \switchedSysAdj \quad \left\{\quad \begin{aligned}
        -\mE_{\switch(t)}^\T\dot{\adState}(t)  &= \mA_{\switch(t)}^\T \adState(t) + \mC_{\switch(t)}^\T (\out(t)-\outVarDes(t)+\multiplier(t)), & \adState(\finalTime) &= \zeroVec, \\
        \mE_{\switch(t^-)}^\T\adState(t^-)&=\mK_{\switch(t^+),\switch(t^-)}^\T\mE_{\switch(t^+)}^\T\adState(t^+),\\
        \adOut(t) &= \mB_{\switch(t)}^\T \adState(t),\\
    \end{aligned}\right.
\end{equation}
which we refer to as the \emph{adjoint system} associated with \eqref{eq:SLS:FOM}.
As in \cite{KarMUV26}, integration by parts on the switching intervals yields the terminal and transition conditions.
The full and reduced optimality systems are stated in the next two propositions.
\begin{proposition}[\MY-regularized \FOM-\OCP optimality system]\label{rem:FOM:opsysoutconsMY}
    The \MY-regularized \FOM-\OCP~\eqref{eq:OPT:FOMoutconsMY} admits a unique solution $(\outOpt,\inpOpt)$. Together with the optimal state $\stateOpt$, optimal adjoint $\adStateOpt$, optimal adjoint output $\adOutOpt$, and the approximate Lagrange multiplier $\multiplierOpt\in L^2(0,\finalTime;\R^\outDim)$ it satisfies
     \begin{subequations}
        \label{eq:OptCond:FOMoutcons}
        \begin{empheq}[left=\empheqlbrace]{align}
            \mE_{\switch(t)}\dstateOpt(t) &= \mA_{\switch(t)}\stateOpt(t) + \mB_{\switch(t)}\inpOpt(t), \quad\stateOpt(0) = \zeroVec, \label{eq:OptCond:FOM:state1outcon}\\
            \stateOpt(t^+) &= \mK_{\switch(t^+),\switch(t^-)}\stateOpt(t^-), \label{eq:OptCond:FOM:state2outcon}\\
            \outOpt(t) &= \mC_{\switch(t)}\stateOpt(t), \label{eq:OptCond:FOM:state3outcon}\\[0.3ex]
            -\mE_{\switch(t)}^\T\dadStateOpt(t) &= \begin{aligned}[t] &\mA_{\switch(t)}^\T\adStateOpt(t) + \mC_{\switch(t)}^\T(\outOpt(t)-\outVarDes(t)+\multiplierOpt(t)),\\ &\adStateOpt(\finalTime) = \zeroVec,\end{aligned} \label{eq:OptCond:FOM:Adstate1outcon}\\
           \mE_{\switch(t^-)}^\T\adStateOpt(t^{-}) &= \mK_{\switch(t^+),\switch(t^-)}^\T\mE_{\switch(t^+)}^\T\adStateOpt(t^+), \label{eq:OptCond:FOMMY:adstatetrans}\\
            \adOutOpt(t) &= \mB^\T_{\switch(t)}\adStateOpt(t), \label{eq:OptCond:FOMAdstate3outcon}\\[0.3ex]
            \multiplierOpt(t) &= \reg\big(\outOpt-\Proj(\outOpt)\big)(t), \label{eq:OptCond:complementary}\\
            0 & \leq \langle \lambda(\inpOpt-\inpVarDes) + \adOutOpt, \mv-\inpOpt\rangle_{L^2(0,\finalTime;\R^\inpDim)} \quad \forall \mv\in\mUad. \label{eq:OptCond:FOM_gradoutcon}
        \end{empheq}
    \end{subequations}
\end{proposition}

\pagebreak[2]
\begin{proposition}[\ROM-\OCP optimality system]\label{rem:ROM:opsysMY}
    %
	The \MY-regularized \ROM-\OCP \eqref{eq:OPT:ROMoutconsMY} admits a unique solution $(\outRedOpt,\inpRedOpt)$. Together with the optimal state $\stateRedOpt$, optimal adjoint $\adStateRedOpt$, optimal adjoint output $\adOutRedOpt$, and the approximate Lagrange multiplier $\multiplierRedOpt\in L^2(0,\finalTime;\R^\outDim)$ it satisfies
    \begin{subequations}
        \label{eq:OptCond:ROMMY}
        \begin{empheq}[left=\empheqlbrace]{align}
            \mERed_{\switch(t)}\dstateRedOpt(t) &= \mARed_{\switch(t)}\stateRedOpt(t) + \mBRed_{\switch(t)}\inpRedOpt(t), \quad \stateRedOpt(0) = \zeroVec,\label{eq:OptCond:ROMMY:state1}\\
            \stateRedOpt(t^+) &= \reduce{\mK}_{\switch(t^+),\switch(t^-)}\stateRedOpt(t^-), \label{eq:OptCond:ROMMY:state2}\\
            \outRedOpt(t) &= \mCRed_{\switch(t)}\stateRedOpt(t), \label{eq:OptCond:ROMMY:state3}\\[0.3ex]
            -\mERed_{\switch(t)}^\T\dadStateRedOpt(t) &= \begin{aligned}[t] &\mARed_{\switch(t)}^\T\adStateRedOpt(t) + \mCRed_{\switch(t)}^\T(\outRedOpt(t)-\outVarDes(t)+\multiplierRedOpt(t)),\\ &\adStateRedOpt(\finalTime) = \zeroVec,\end{aligned}\label{eq:OptCond:ROMMY:Adstate1}\\
            \mERed_{\switch(t^-)}^\T\adStateRedOpt(t^{-}) &= \reduce{\mK}_{\switch(t^+),\switch(t^-)}^\T\mERed_{\switch(t^+)}^\T\adStateRedOpt(t^+), \label{eq:OptCond:ROMMY:adstatetrans}\\
            \adOutRedOpt(t) &= \mBRed^\T_{\switch(t)}\adStateRedOpt(t), \label{eq:OptCond:ROMMY:Adstate3}\\[0.3ex]
            \multiplierRedOpt(t) &= \reg\big(\outRedOpt-\Proj(\outRedOpt)\big)(t), \label{eq:OptCond:ROMMY:complementary}\\
            0 &\leq \langle \lambda(\inpRedOpt-\inpVarDes) + \adOutRedOpt, \mv-\inpRedOpt\rangle_{L^2(0,\finalTime;\R^\inpDim)} \quad \forall \mv\in\mUad. \label{eq:OptCond:ROMMY_grad}
        \end{empheq}
    \end{subequations}
\end{proposition}

To establish an output-based a posteriori error estimator, we first collect properties of the regularized multiplier.
\begin{proposition}[Properties of the approximate Lagrange multiplier]\label{rem:multiplier:properties}
    Let $\out\in L^2(0,\finalTime;\R^\outDim)$ and define $\multiplier(\out)=\reg(\out-\Proj(\out))$.
    Then the following properties hold:
    \begin{enumerate}
        \item[\textup{(i)}]$ \langle \multiplier(\out), \mv-\Proj(\out) \rangle_{L^2} \leq 0$ for all $\mv\in \mYad$.
        \item[\textup{(ii)}]$ \langle \multiplier(\out)-\multiplier(\hat\out), \Proj(\hat\out)-\Proj(\out) \rangle_{L^2} \leq 0$ for all $\out, \hat \out \in L^2(0,T;\R^\outDim)$.
    \end{enumerate}
\end{proposition}
\begin{proof}
    Property \textup{(i)} follows immediately from the characterization of the metric projection and the definition of $\multiplier$.
    To prove \textup{(ii)}, apply this characterization first to $\out$ with $\mv=\Proj(\hat\out)$ and then to $\hat\out$ with $\mv=\Proj(\out)$.
    Combining the resulting inequalities using $\reg>0$ yields \textup{(ii)}.
\end{proof}

Now, we estimate the error in the optimal control and output.
\begin{theorem}[Output-based a posteriori error estimate for $\inpOpt,\outOpt$]\label{lemma:apost_control}
    \leavevmode\\
    Let $(\outOpt,\inpOpt, \multiplierOpt)$ and $(\outRedOpt, \inpRedOpt, \multiplierRedOpt)$ be the optimal output, input, and approximate Lagrange multiplier of the \FOM-\OCP \eqref{eq:OPT:FOMoutconsMY} and the \ROM-\OCP \eqref{eq:OPT:ROMoutconsMY}, respectively. Define the error 
    \begin{align}
    	\label{eq:err:OCP:ROM}
        e_{\outOpt,\inpOpt} &\vcentcolon= \bigg(\|\outRedOpt-\outOpt \|^2_{L^2} + \lambda\|\inpRedOpt-\inpOpt \|^2_{L^2}+\tfrac{1}{\reg}\|\multiplierRedOpt-\multiplierOpt \|^2_{L^2}\bigg)^{\nicefrac{1}{2}}.
	\end{align}
	Then $e_{\outOpt,\inpOpt} \leq \Delta_{\outOpt,\inpOpt}(\outRedOpt,\inpRedOpt,\multiplierRedOpt)$ with 
	\begin{align}
		\label{eq:abstract_control_bound}
        \Delta_{\outOpt,\inpOpt}(\outRedOpt,\inpRedOpt,\multiplierRedOpt)\;&\vcentcolon=\; \bigg(\tfrac{1}{\lambda}\|\adOut(\outRedOpt,\multiplierRedOpt)-\adOutRedOpt\|^2_{L^2}+(1+\reg)\|\out(\inpRedOpt)-\outRedOpt\|^2_{L^2}\bigg)^{\nicefrac{1}{2}}.
    \end{align}
    %
\end{theorem}
\begin{proof}
    We adapt the proof from \cite{KarMUV26} to \SLS with state-transition maps, Petrov--Galerkin \ROMs, and \MY-regularized output constraints.
    Choose $\mv=\inpRedOpt$ in \eqref{eq:OptCond:FOM_gradoutcon} and $\mv=\inpOpt$ in \eqref{eq:OptCond:ROMMY_grad}.
    Adding the resulting inequalities gives
    \begin{equation*}
        \lambda\|\inpRedOpt-\inpOpt\|^2_{L^2}\;\leq\; \langle \adOutRedOpt-\adOutOpt,\inpOpt-\inpRedOpt \rangle_{L^2}.
    \end{equation*}
    Adding and subtracting $\adOut(\outRedOpt, \multiplierRedOpt)$, which is the output of the \FOM adjoint system \eqref{eq:SLS:FOM_Adj} with inputs $\outRedOpt, \multiplierRedOpt$, leads to
    \begin{equation}\label{eq:ocp_est_helper}
        \lambda\|\inpRedOpt-\inpOpt\|^2_{L^2}\leq \langle \adOutRedOpt-\adOut(\outRedOpt,\multiplierRedOpt),\inpOpt-\inpRedOpt \rangle_{L^2}+\langle\adOut(\outRedOpt,\multiplierRedOpt)-\adOutOpt,\inpOpt-\inpRedOpt \rangle_{L^2}.
    \end{equation}
    By the definition of the adjoint output, the second term on the right-hand side is given as
    \begin{align*}
        \langle\adOut(\outRedOpt,\multiplierRedOpt)-\adOutOpt,\inpOpt-\inpRedOpt \rangle_{L^2}
        =&   \int_{0}^\finalTime \langle\adState(\outRedOpt,\multiplierRedOpt)- \adStateOpt,\mB_{\switch}(\inpOpt-\inpRedOpt )\rangle_{\R^\stateDim}(t)\ \dt
    \end{align*}
    Using the \FOM \eqref{eq:SLS:FOM} with input $\inpOpt-\inpRedOpt$ and integrating by parts on each switching interval leads to
    \begin{align}\nonumber
        &\langle\adOut(\outRedOpt,\multiplierRedOpt)-\adOutOpt,\inpOpt-\inpRedOpt \rangle_{L^2}\\ \nonumber
        &=\int_{0}^\finalTime \bigl\langle\adState(\outRedOpt,\multiplierRedOpt)-\adStateOpt,
        \mE_{\switch}\ddt(\stateOpt-\state(\inpRedOpt))-\mA_{\switch}(\stateOpt-\state(\inpRedOpt))\bigr\rangle_{\R^\stateDim}(t)\ \dt\\ \label{eq:est_arg}
        &=\int_{0}^\finalTime \bigl\langle-\mE_{\switch}^\T\ddt(\adState(\outRedOpt,\multiplierRedOpt)-\adStateOpt)
        -\mA_{\switch}^\T(\adState(\outRedOpt,\multiplierRedOpt)-\adStateOpt),\stateOpt-\state(\inpRedOpt)\bigr\rangle_{\R^\stateDim}(t)\ \dt\\ \nonumber
        &\quad+\sum_{k=0}^{\nrSwitches}\left[
        \bigl\langle\mE_{\switch}^\T(\adState(\outRedOpt,\multiplierRedOpt)-\adStateOpt),
        \stateOpt-\state(\inpRedOpt)\bigr\rangle_{\R^\stateDim}(t)
        \right]_{t=t_k^+}^{t=t_{k+1}^-}.
    \end{align}
    The final-time term vanishes by the terminal conditions in \eqref{eq:SLS:FOM_Adj} and \eqref{eq:OptCond:FOM:Adstate1outcon}, and the initial-time term vanishes since $(\stateOpt-\state(\inpRedOpt))(0^+)=\zeroVec$ (see \eqref{eq:SLS:FOM}).
    The switching terms cancel as both full-order adjoint states satisfy the transition conditions in \eqref{eq:SLS:FOM_Adj} and \eqref{eq:OptCond:FOMMY:adstatetrans}, respectively.
    Applying the \FOM adjoint equation \eqref{eq:SLS:FOM_Adj} to $\adState(\outRedOpt,\multiplierRedOpt)-\adStateOpt$ yields
    \begin{align}\label{eq:est_arg2}
        \begin{aligned}
        \langle\adOut(\outRedOpt,\multiplierRedOpt)-\adOutOpt,&\inpOpt-\inpRedOpt \rangle_{L^2}
        \\
        =&\;\langle\outRedOpt-\outOpt,\outOpt-\out(\inpRedOpt)\rangle_{L^2}+\langle\multiplierRedOpt-\multiplierOpt,\outOpt-\out(\inpRedOpt)\rangle_{L^2}.
        \end{aligned}
    \end{align}
   The first term is estimated via the Cauchy--Schwarz and Young inequalities, for $\varepsilon_1>0$, as
    \begin{align}%
        \nonumber
        \langle\outRedOpt-\outOpt,\outOpt-\out(\inpRedOpt)\rangle_{L^2}\;=&\;
        -\|\outOpt-\outRedOpt\|^2_{L^2}+\langle\outRedOpt-\outOpt,\outRedOpt-\out(\inpRedOpt)\rangle_{L^2}\\ \nonumber
        \leq &\; (\tfrac{\varepsilon_1}{2}-1)\|\outOpt-\outRedOpt\|^2_{L^2}+\tfrac{1}{2\varepsilon_1}\|\outRedOpt-\out(\inpRedOpt)\|^2_{L^2}.
    \end{align}
    Using \Cref{rem:multiplier:properties}\textup{(ii)} and Young's inequality, the second term is estimated for $\varepsilon_2>0$ as
    \begin{align*}
        \langle\multiplierRedOpt-\multiplierOpt,\outOpt-\out(\inpRedOpt)\rangle_{L^2}\;&=\;
        \langle\multiplierRedOpt-\multiplierOpt,\outOpt-\outRedOpt\rangle_{L^2}+
        \langle\multiplierRedOpt-\multiplierOpt,\outRedOpt-\out(\inpRedOpt)\rangle_{L^2}\\
        &=\;-\tfrac{1}{\reg}\|\multiplierRedOpt-\multiplierOpt\|_{L^2}^2+\langle\multiplierRedOpt-\multiplierOpt,\Proj(\outOpt)-\Proj(\outRedOpt)\rangle_{L^2}\\
        &\phantom{=\;}\;+\langle\multiplierRedOpt-\multiplierOpt,\outRedOpt-\out(\inpRedOpt)\rangle_{L^2}\\
        &\leq\;  \left(\tfrac{\varepsilon_2}{2} -\tfrac{1}{\reg}\right)\|\multiplierRedOpt-\multiplierOpt\|_{L^2}^2+\tfrac{1}{2\varepsilon_2}\|\outRedOpt-\out(\inpRedOpt)\|^2_{L^2}.
    \end{align*}
    Combining both estimates leads to
    \begin{align*}
        \langle\adOut(\outRedOpt,\multiplierRedOpt)-\adOutOpt,\inpOpt-\inpRedOpt \rangle_{L^2}\; & \leq\; \left(\tfrac{1}{2\varepsilon_1}+\tfrac{1}{2\varepsilon_2}\right)\|\outRedOpt-\out(\inpRedOpt)\|^2_{L^2}+ \left(\tfrac{\varepsilon_2}{2}-\tfrac{1}{\reg}\right)\|\multiplierRedOpt-\multiplierOpt\|_{L^2}^2\\
        &+\; \left(\tfrac{\varepsilon_1}{2}-1\right)\|\outOpt-\outRedOpt\|^2_{L^2}.
    \end{align*}
    For the first term on the right-hand side in \eqref{eq:ocp_est_helper} we use the weighted Young's inequality with $\varepsilon_3>0$ to obtain
    \begin{align*}
        \langle\adOutRedOpt-\adOut(\outRedOpt,\multiplierRedOpt),\inpOpt-\inpRedOpt \rangle_{L^2}\;\leq\; \tfrac{1}{2\varepsilon_3}\| \adOutRedOpt-\adOut(\outRedOpt,\multiplierRedOpt)\|^2_{L^2}+\tfrac{\varepsilon_3}{2}\|\inpOpt-\inpRedOpt\|^2_{L^2}.
    \end{align*}
    Combining all inequalities yields the estimate
    \begin{align}
        \label{eq:ocp_est_helper3}
        \begin{aligned}
            &\left(1-\tfrac{\varepsilon_1}{2}\right)\|\outOpt-\outRedOpt\|^2_{L^2}+\left( \tfrac{1}{\reg}-\tfrac{\varepsilon_2}{2}\right)\|\multiplierRedOpt-\multiplierOpt\|_{L^2}^2+\left(\lambda-\tfrac{\varepsilon_3}{2}\right)\|\inpRedOpt-\inpOpt\|^2_{L^2}\\
            &\leq\; \left(\tfrac{1}{2\varepsilon_1}+\tfrac{1}{2\varepsilon_2}\right)\|\outRedOpt-\out(\inpRedOpt)\|^2_{L^2} +\tfrac{1}{2\varepsilon_3}\| \adOutRedOpt-\adOut(\outRedOpt,\multiplierRedOpt)\|^2_{L^2}.
        \end{aligned}
    \end{align}
    Choosing $\varepsilon_1=1$, $\varepsilon_2=\nicefrac{1}{\reg}$, and $\varepsilon_3=\lambda$, multiplying by two, and taking square roots proves the claim.
\end{proof}

Next, we derive an error estimate for the optimal value function.
\begin{theorem}[A posteriori error estimate for the value function]\label{lemma:apost_valuefun}
    In the situation of \Cref{lemma:apost_control}, we have
    \begin{equation}\label{eq:value:convex:interval}
        \eta_-\leq
        \calJ_\reg(\outOpt,\inpOpt)-\calJ_\reg(\outRedOpt,\inpRedOpt)
        \leq\eta_+,
    \end{equation}
    where
    \begin{equation}\label{eq:value:eta}
        \eta_\pm\vcentcolon=
        \bigl\langle\outRedOpt-\outVarDes+\multiplierRedOpt,
        \out(\inpRedOpt)-\outRedOpt\bigr\rangle_{L^2}
        +\begin{cases}
            -\tfrac{1}{2\lambda}\|\adOut(\outRedOpt,\multiplierRedOpt)-\adOutRedOpt\|_{L^2}^2,&\text{for }\eta_-,\\
            \tfrac{1+\reg}{2}\|\out(\inpRedOpt)-\outRedOpt\|_{L^2}^2,&\text{for }\eta_+.
        \end{cases}
    \end{equation}
    Consequently, the value function error  
    \begin{equation}\label{eq: val est 0}
        e_{\calJ_\reg}\vcentcolon=
        |\calJ_\reg(\outOpt,\inpOpt)-\calJ_\reg(\outRedOpt,\inpRedOpt)|
    \end{equation}
    satisfies
    \begin{equation}\label{eq:value:convex:bound}
        e_{\calJ_\reg}\leq
        \Delta_{\calJ_\reg}(\outRedOpt,\inpRedOpt,\multiplierRedOpt)
        \vcentcolon=\max\{|\eta_-|,|\eta_+|\}.
    \end{equation}
\end{theorem}
\begin{proof}
    For fixed $\inp$, the output-dependent part of $\calJ_\reg$ is convex,
    with gradient $\out-\outVarDes+\multiplier(\out)$.
    This gradient is $(1+\reg)$-Lipschitz, since
    $\out\mapsto\out-\Proj(\out)$ is nonexpansive.
    Turning to the lower bound, for every $\mv\in\mUad$, convexity in the output and the quadratic control cost give
    \begin{align*}
        \calJ_\reg(\out(\mv),\mv)-\calJ_\reg(\outRedOpt,\inpRedOpt)
        &\geq \bigl\langle\outRedOpt-\outVarDes+\multiplierRedOpt,
        \out(\mv)-\outRedOpt\bigr\rangle_{L^2}\\
        &\quad +\lambda\langle\inpRedOpt-\inpVarDes,\mv-\inpRedOpt\rangle_{L^2}
        +\tfrac{\lambda}{2}\|\mv-\inpRedOpt\|_{L^2}^2.
    \end{align*}
    Split $\out(\mv)-\outRedOpt$ into
    $\out(\mv)-\out(\inpRedOpt)$ and $\out(\inpRedOpt)-\outRedOpt$.
    Linearity of the \FOM input-to-output map, the adjoint identity
    as in \eqref{eq:est_arg2}, and the reduced optimality condition
    \eqref{eq:OptCond:ROMMY_grad} then imply
    \begin{align*}
        \calJ_\reg(\out(\mv),\mv) - \calJ_\reg(\outRedOpt,\inpRedOpt) &\geq
        \bigl\langle\outRedOpt-\outVarDes+\multiplierRedOpt,
        \out(\inpRedOpt)-\outRedOpt\bigr\rangle_{L^2}\\
        &\qquad+
        \bigl\langle\adOut(\outRedOpt,\multiplierRedOpt)-\adOutRedOpt,
        \mv-\inpRedOpt\bigr\rangle_{L^2}
        +\tfrac{\lambda}{2}\|\mv-\inpRedOpt\|_{L^2}^2\\
        &\geq
        \bigl\langle\outRedOpt-\outVarDes+\multiplierRedOpt,
        \out(\inpRedOpt)-\outRedOpt\bigr\rangle_{L^2}
        \\&\qquad-\tfrac{1}{2\lambda}
        \|\adOut(\outRedOpt,\multiplierRedOpt)-\adOutRedOpt\|_{L^2}^2
        =\eta_-,
    \end{align*}
    %
    where the last step is Young's inequality.
    Taking $\mv=\inpOpt$ proves the lower bound
    in \eqref{eq:value:convex:interval}.
    For the upper bound,
    optimality of $\inpOpt$, feasibility of $\inpRedOpt\in\mUad$,
    and the Lipschitz continuity of the output gradient give
    \begin{align*}
        \calJ_\reg(\outOpt,\inpOpt)-\calJ_\reg(\outRedOpt,\inpRedOpt) &\leq
        \calJ_\reg(\out(\inpRedOpt),\inpRedOpt)
        -\calJ_\reg(\outRedOpt,\inpRedOpt)\\
        &\leq
        \bigl\langle\outRedOpt-\outVarDes+\multiplierRedOpt,
        \out(\inpRedOpt)-\outRedOpt\bigr\rangle_{L^2}
        +\tfrac{1+\reg}{2}\|\out(\inpRedOpt)-\outRedOpt\|_{L^2}^2
        =\eta_+.
    \end{align*}
    This proves the upper bound in \eqref{eq:value:convex:interval}.
    Taking the maximum of the absolute values of the two endpoints
    yields \eqref{eq:value:convex:bound}.
\end{proof}

To estimate how the reduced control performs in the \FOM cost, note that the map $\inp\mapsto\calJ_\reg(\out(\inp),\inp)$ is $\lambda$-strongly convex, with gradient $\lambda(\inp-\inpVarDes)+\adOut(\out(\inp),\multiplier(\out(\inp)))$. Thus, using the reduced optimality condition \eqref{eq:OptCond:ROMMY_grad} with $\mv=\inpOpt$ yields
\begin{align*}
	0&\leq\calJ_\reg(\out(\inpRedOpt),\inpRedOpt)-\calJ_\reg(\outOpt,\inpOpt)\\
	&\leq\bigl\langle\lambda(\inpRedOpt-\inpVarDes) +\adOut(\out(\inpRedOpt),\multiplier(\out(\inpRedOpt))), \inpRedOpt-\inpOpt\bigr\rangle_{L^2} - \tfrac{\lambda}{2}\|\inpRedOpt-\inpOpt\|_{L^2}^2\\
	&\leq\bigl\langle\adOut(\out(\inpRedOpt),\multiplier(\out(\inpRedOpt))) -\adOutRedOpt,\inpRedOpt-\inpOpt\bigr\rangle_{L^2} - \tfrac{\lambda}{2}\|\inpRedOpt-\inpOpt\|_{L^2}^2.
\end{align*}
The bounds in \Cref{lemma:apost_control} and \Cref{lemma:apost_valuefun} can be evaluated without computing the \FOM optimal control.
Their direct evaluation requires one forward \FOM simulation driven by $\inpRedOpt$, to compute $\out(\inpRedOpt)$, and one backward \FOM adjoint simulation driven by, e.g., $(\outRedOpt,\multiplierRedOpt)$, to compute $\adOut(\outRedOpt,\multiplierRedOpt)$.
This direct evaluation uses the full-order forward and adjoint models without requiring full-order optimization, which motivates the term full-order-evaluation bounds. The second numerical experiment compares these bounds with an online-efficient \PBR approximation.
\begin{remark}[Augmented-Lagrangian regularization]\label{rem:AL:shifts}
    Our error estimation framework can also be adapted to augmented Lagrangian methods (see, e.g., \cite{HPUU09}). In this case, for a fixed shift $\shift\in L^2(0,\finalTime;\R^\outDim)$, the \MY penalty in \eqref{eq:OPT:FOMoutconsMY} is replaced by
    $\tfrac{\reg}{2}\operatorname{dist}_{L^2}^2(\out+\shift/\reg,\mYad)
    -\tfrac{1}{2\reg}\|\shift\|_{L^2}^2$.
    The effective multiplier is
    $\multiplier_{\reg,\shift}(\out)=\reg(\out+\shift/\reg-\Proj(\out+\shift/\reg))$.
    Then, the estimates in
    \Cref{lemma:apost_control,lemma:apost_valuefun} extend unchanged, using
    this multiplier in the adjoint input.
\end{remark}
%

\section{Piecewise balanced reduction for the optimal control problem}
\label{sec:BT_for_OCP}
%
The analysis from \Cref{sec:outputbased_error_est} can be applied to any Petrov--Galerkin projection-based \ROM \eqref{eq:SLS:ROM} that approximates the forward and adjoint input-output maps.
We show in \Cref{subsec:duality} that the \PBR from~\cite{ManU26} is a well-suited tool for reducing the \OCP, since the duality between observability and controllability makes \PBR commute with the derivation of the adjoint system: reducing the \FOM first yields the same adjoint system as deriving the full-order adjoint first and then applying \PBR.

Throughout this section, we impose the following assumption.
\begin{assumption}\label{ass:BT}
    For each $j=1,\ldots,L$, all eigenvalues of $\mE_j^{-1}\mA_j$ lie in the open left half-plane.
\end{assumption}
\Cref{ass:BT} ensures that each mode of the \SLS is asymptotically stable.
Although this condition alone does not guarantee asymptotic stability for the \SLS under arbitrary switching signals, it remains useful for introducing the notion of generalized Gramians defined through Lyapunov-type equations, which are fundamental to the \PBR method, which we recall in \Cref{subsec:PBR}. 

\subsection{The \PBR for \SLS}\label{subsec:PBR}
Given the \SLS~\eqref{eq:SLS:FOM}, we consider for each $i\in\switchingSet$ the pair of \GLEs
\begin{subequations}
    \label{eqn:GLE:b}
    \begin{align}
        \mA_i \fcP_{i,\state}\mE_i^\T+\mE_i\fcP_{i,\state}\mA_i^\T+\sum_{j=1}^{L}\left(\mN_{i,j,\state} \fcP_{i,\state}\mN_{i,j,\state}^\T+\mB_{i,j,\state}\mB_{i,j,\state}^\T\right) &= \zeroVec,\label{eqn:GLE:reach:b}\\
        \mA_i^\T \fcQ_{i,\state}\mE_i + \mE_i^\T\fcQ_{i,\state}\mA_i+\sum_{j=1}^{L}\left(\mN^\T_{i,j,\state} \fcQ_{i,\state}\mN_{i,j,\state} + \mC_{i,j,\state}^\T\mC_{i,j,\state}\right) &= \zeroVec,\label{eqn:GLE:observ:b}
    \end{align}
\end{subequations}
where for all $i,j\in\switchingSet$
\begin{align*}
	\mB_{i,j,\state} &\vcentcolon= \mE_i\mK_{i,j}\mE_j^{-1}\mB_j, &
	\mC_{i,j,\state} &\vcentcolon= \mC_j\mK_{j,i}, &
	\mN_{i,j,\state} &\vcentcolon= \mE_i\mK_{i,j}\mE_j^{-1}\mA_j\mK_{j,i}-\mA_i.
\end{align*}
Note that we use the subscript $\state$ to emphasize that the matrices are computed for the forward system~\eqref{eq:SLS:FOM}. In contrast, if we compute the matrices for the adjoint system, we write $\fcP_{i,\adState}$ and similarly for the other matrices.
The existence and uniqueness of the solutions to~\eqref{eqn:GLE:reach:b}, and analogously to \eqref{eqn:GLE:observ:b}, can be studied using Kronecker algebra.
To this end, define $\lyapOper_i, \gleOper_{i,\state}\in\R^{\stateDim^2\times\stateDim^2}$ via
\begin{align}\label{eqn:GLE:oper}
    \lyapOper_i \vcentcolon= \mE_i\otimes\mA_i +\mA_i\otimes\mE_i \qquad\text{and}\qquad
    \gleOper_{i,\state} \vcentcolon= \sum_{j=1}^{L}\mN_{i,j,\state}\otimes\mN_{i,j,\state}.
\end{align}
The vectorized form of~\eqref{eqn:GLE:reach:b} is then $(\lyapOper_i+\gleOper_{i,\state})\vec2(\fcP_{i,\state})=-\vec2(\sum_{j=1}^{L}\mB_{i,j,\state}\mB_{i,j,\state}^\T)$.
Hence, \eqref{eqn:GLE:reach:b} is uniquely solvable if and only if $\lyapOper_i+\gleOper_{i,\state}$ is nonsingular.
\begin{assumption}\label{ass:uniq:GLE}
    The \GLEs in \eqref{eqn:GLE:b} are uniquely solvable for all $i\in\switchingSet$.
\end{assumption}
The solutions $\fcP_{i,\state}$ and $\fcQ_{i,\state}$ are symmetric positive semidefinite due to \Cref{ass:BT} and satisfy $\calR = \text{range}(\fcP_{i,\state})$ and $\calO = \text{range}(\mE_i^\T\fcQ_{i,\state}\mE_i)$ for $i=1,\ldots,L$, where $\calR$ and $\calO$ are the reachable and observable sets of the \SLS, respectively (see \cite[Sec.~2.1]{ManU26} and \cite[Thm.~3]{PonGB20}).
We therefore refer to the $L$ pairs $(\fcP_{i,\state},\fcQ_{i,\state})$ as the reachability and observability matrices of the \SLS.

Let $\mS_{i,\state}$ and $\mR_{i,\state}$ be Cholesky factors satisfying $\fcP_{i,\state}=\mS_{i,\state}\mS_{i,\state}^\T$ and $\fcQ_{i,\state}=\mR_{i,\state}\mR_{i,\state}^\T$.
Since the observability matrix of the normalized subsystem is $\mE_i^\T\fcQ_{i,\state}\mE_i$ with associated Cholesky factor $\mE_i^\T\mR_{i,\state}$, we compute the mass-weighted singular value decomposition
\begin{equation}\label{eqn:SVD}
    \mS_{i,\state}^\T\mE_i^\T\mR_{i,\state} = \left[\fU_{i,\state,1}, \fU_{i,\state,2}\right]\begin{bmatrix}
        \mathbf{\Sigma}_{i,\state,1}&\zeroVec\\
        \zeroVec&\mathbf{\Sigma}_{i,\state,2}
    \end{bmatrix}[\mPsi_{i,\state,1},\mPsi_{i,\state,2}]^\T,
\end{equation}
where the first block contains the $r$ retained positive singular values.
The projection matrices~$\mV_{i,\state}$ and $\mW_{i,\state}$ are then given by 
\begin{equation}\label{eqn:PBR:proj}
    \mV_{i,\state} = \mS_{i,\state}\fU_{i,\state, 1}\mathbf{\Sigma}^{-1/2}_{i,\state, 1} \quad\text{and}\quad
    \mW_{i,\state} = \mR_{i,\state}\mPsi_{i,\state,1}\mathbf{\Sigma}^{-1/2}_{i,\state,1}.
\end{equation}
These matrices satisfy
\begin{equation*}
    \mW_{i,\state}^\T\mE_i\mV_{i,\state}
    =\mathbf{\Sigma}_{i,\state,1}^{-1/2}\mPsi_{i,\state,1}^\T
    (\mR_{i,\state}^\T\mE_i\mS_{i,\state})
    \fU_{i,\state,1}\mathbf{\Sigma}_{i,\state,1}^{-1/2}
    =\mI_r.
\end{equation*}
Thus, the reduced descriptor matrices satisfy $\reduce{\mE}_i=\mI_r$, even when $\mE_i\neq\mI_\stateDim$.
Applying these projection matrices in~\eqref{eqn:SLS:ROM:matrices} yields the \ROM introduced in~\eqref{eq:SLS:ROM} with transition matrices $\reduce{\mK}_{i,j}=\mW_{i,\state}^\T\mE_i\mK_{i,j}\mV_{j,\state}$. This procedure is called \PBR for \SLS.
Its stability properties are discussed in \cite[Sec.~4.2]{ManU26}.
Moreover, \PBR provides an a posteriori estimator for the $L^2$-norm of the output error \cite[Thm.~4.3]{ManU26}:
\begin{align}\label{eqn:PBR:bound}
    \begin{aligned}
        \|\out(\inp)-\outRed(\inp)\|_{L^2(0,T;\R^{\outDim})}\;\le&\;\Delta(\mathbf{\Sigma}_{\state,2},\stateRed,\lambda_1,\inp),\quad\text{where}\\
        \Delta(\mathbf{\Sigma}_{\state,2},\stateRed,\lambda_1,\inp)\;\vcentcolon=\;\sum_{i=\stateDimRed+1}^{\stateDim}\Big(&2\|\inp\|_{L^2}f_{1,i}(\sigma_{1},\ldots,\sigma_{L})
        \\
        +&\;f_{2,i}(\stateRed_{-},\stateRed_{+})+f_{3,i}(\lambda_1(\mM_c),\lambda_1(\mM_o))\Big),
    \end{aligned}
\end{align}
and $\Delta(\mathbf{\Sigma}_{\state,2},\stateRed,\lambda_1,\inp)$ depends on the sum of three main contributions.
The terms $f_{1,i}$ account for the neglected singular values in~\eqref{eqn:SVD}.
Through $f_{2,i}$, the effect of discontinuities in the reduced state at the switching instants is quantified.
Finally, the terms involving $f_{3,i}$ account for violations of the relevant linear matrix inequalities as measured by the largest positive eigenvalue. The derivation of the bound, the precise definitions of all quantities involved, and a detailed discussion of its properties can be found in~\cite[Sec.~4]{ManU26}. We note that the error estimator $\Delta(\mathbf{\Sigma}_{\state,2},\stateRed,\lambda_1,\inp)$ in \eqref{eqn:PBR:bound} is not computationally efficient to evaluate, as it requires an a posteriori computation of trajectories in the full-order space of dimension $\stateDim$. In practical implementations, the summation up to $\stateDim$ is therefore replaced by a truncated sum up to $\tilde\stateDimRed$, where $\stateDimRed<\tilde\stateDimRed\ll\stateDim$. Provided that the series terms decay sufficiently rapidly to zero, or converge to a limiting value (noting that some contributions may stagnate as the dimension increases), the resulting quantity serves as a reliable error indicator while remaining computationally tractable. 

The estimator simplifies for certain problem classes, for example when the matrices $\mE_i^{-1}\mA_i$ are well-conditioned and have negative-definite Hermitian parts (see \cite[Sec.~3.4]{ManU26}).
In particular, if the \SLS has identity state transitions and admits common balancing bases whose Gramians satisfy the required \LMIs for all modes, the switching and \LMI-defect corrections can be omitted. In this case, we recover the classical \BT error bound
\begin{equation}\label{eqn:PBR:special:bound}
    \|\out(\inp)-\outRed(\inp)\|_{L^2}
    \leq c_r\|\inp\|_{L^2},\qquad
    c_r\vcentcolon=2\operatorname{tr}(\mathbf{\Sigma}_{1,\state,2}).
\end{equation}
Since $c_r$ can be precomputed offline, this bound can be evaluated online without trajectory-dependent correction terms; see~\cite[Thm.~2.4]{ManU26}.
In the general case, a rigorous evaluation of the estimator can be computationally demanding, and the resulting bound may substantially overestimate the true reduction error; consequently, it is often preferable to trade strict rigor for improved computational efficiency.

\subsection{Duality of the \PBR for the optimal control problem}\label{subsec:duality}
\begin{figure}[ht]
	\centering
	\includegraphics{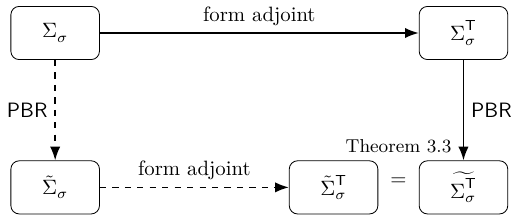}
	\caption{First-optimize-then-reduce (solid path) vs.~first-reduce-then-optimize (dashed path).}
	\label{fig:optimize-vs-reduce}
\end{figure}

%
So far, we obtained the optimality system in \Cref{rem:ROM:opsysMY} and the resulting error estimators presented in  \Cref{lemma:apost_control,lemma:apost_valuefun} by the \emph{first-reduce-then-optimize} approach, illustrated by the dashed arrows in~\Cref{fig:optimize-vs-reduce}. In this approach, we first compute the \ROM~$\switchedSysRed$ and then derive its optimality system, which includes the adjoint system~$\switchedSysRed^{\T}$ associated with the \ROM.
Alternatively, we first derive the optimality system for the \FOM, yielding the adjoint system~$\switchedSysAdj$, and then apply \MOR to the adjoint system to obtain the reduced adjoint system~$\widetilde{\switchedSysAdj}$. This approach, called \emph{first-optimize-then-reduce}, is illustrated by the solid arrows in~\Cref{fig:optimize-vs-reduce}.
A natural question is whether these two approaches commute, that is, whether
$\widetilde{\switchedSysAdj}\equiv\switchedSysRed^{\T}$.
The following theorem establishes that this is indeed the case for~\PBR.
\begin{theorem}[Commutativity of the \OCP-\PBR and \PBR-\OCP]\label{teo:duality}
    Suppose that \Cref{ass:BT,ass:uniq:GLE} hold. 
    \begin{enumerate}
    	\item\label{item1} Let $(\fcP_{i,\state},\fcQ_{i,\state})_{i=1}^{L}$ and $(\fcP_{i,\adState},\fcQ_{i,\adState})_{i=1}^L$ denote the reachability and observability matrices of~\eqref{eq:SLS:FOM} and~\eqref{eq:SLS:FOM_Adj}, respectively. Then 
    		\begin{equation}
    			\label{eqn:statement}
            	\fcP_{i,\adState} = \fcQ_{i,\state}\qquad\text{and}\qquad \fcQ_{i,\adState} =\fcP_{i,\state} \qquad\text{for }i=1,\ldots,L.
        	\end{equation}
        \item\label{item2} Let $\switchedSysRed^{\T}$ denote the adjoint system corresponding to the \PBR-\ROM $\switchedSysRed$ and let $\widetilde{\switchedSysAdj}$ denote the \PBR-\ROM for the adjoint system~$\switchedSysAdj$. Then $\widetilde{\switchedSysAdj}\equiv\switchedSysRed^{\T}$.
    \end{enumerate}
\end{theorem}
\begin{proof}
    After reversing time in the adjoint system \eqref{eq:SLS:FOM_Adj}, its $L$ pairs of \GLEs read
    \begin{subequations}
        \label{eqn:GLE:b:ad}
        \begin{align}
            \mA^{\T}_i \fcP_{i,\adState}\mE_i+\mE_i^\T\fcP_{i,\adState}\mA_i+\sum_{j=1}^{L}\left(\mN_{i,j,\adState} \fcP_{i,\adState}\mN_{i,j,\adState}^\T+\mC_{i,j,\adState}^{\T}\mC_{i,j,\adState}\right) &= \zeroVec,\label{eqn:GLE:reach:b:ad}\\
            \mA_i \fcQ_{i,\adState}\mE_i^\T + \mE_i\fcQ_{i,\adState}\mA^\T_i+\sum_{j=1}^{L}\left(\mN^\T_{i,j,\adState} \fcQ_{i,\adState}\mN_{i,j,\adState} + \mB_{i,j,\adState}\mB_{i,j,\adState}^{\T}\right) &= \zeroVec,\label{eqn:GLE:observ:b:ad}
        \end{align}
    \end{subequations}
    where $\mK^{\mathrm{ad}}_{j,i}\vcentcolon=\mE_j^{-\T}\mK_{i,j}^{\T}\mE_i^\T$ denotes the full-order adjoint transition from mode $i$ to mode $j$ in reversed time, and for all $i,j\in\switchingSet$
    \begin{align}\label{eqn:adj:mat}
        \begin{aligned}
            \mB_{i,j,\adState} &\vcentcolon= (\mK^{\mathrm{ad}}_{j,i})^\T\mB_j, \\
            \mC_{i,j,\adState} &\vcentcolon= \mC_j\mE_j^{-1}(\mK^{\mathrm{ad}}_{i,j})^\T\mE_i, \\
			\mN_{i,j,\adState} &\vcentcolon= \mE_i^\T\mK^{\mathrm{ad}}_{i,j}\mE_j^{-\T}\mA_j^\T\mK^{\mathrm{ad}}_{j,i}-\mA_i^\T.
        \end{aligned}
    \end{align}
    Substituting the adjoint transition matrices $\mK^{\mathrm{ad}}_{j,i}$ into \eqref{eqn:adj:mat}, we find
    \begin{align}\label{eqn:equ}
        \begin{gathered}
            \mB_{i,j,\adState} = \mE_i\mK_{i,j}\mE_j^{-1}\mB_j = \mB_{i,j,\state},\qquad 
            \mC_{i,j,\adState} = \mC_j\mE_j^{-1}\mE_j\mK_{j,i}\mE_i^{-1}\mE_i = \mC_{i,j,\state},\\
            \mN_{i,j,\adState} = \mK_{j,i}^{\T}\mA_j^{\T}\mE_j^{-\T}\mK_{i,j}^{\T}\mE_i^\T-\mA_i^{\T} = \big(\mE_i\mK_{i,j}\mE_j^{-1}\mA_j\mK_{j,i}-\mA_i\big)^{\T} = \mN_{i,j,\state}^{\T}.
        \end{gathered}
    \end{align}
    Substituting the identities in \eqref{eqn:equ} into \eqref{eqn:GLE:b:ad} and invoking the uniqueness assumption (\Cref{ass:uniq:GLE}) yields \eqref{eqn:statement}. 
    
    For the second statement, we observe that \ref{item1} allows us to choose $\mS_{i,\adState}=\mR_{i,\state}$ and $\mR_{i,\adState}=\mS_{i,\state}$.
    Since the adjoint mass matrix is $\smash{\mE_i^\T}$, its mass-weighted factor product satisfies
    \begin{equation*}
        \mS_{i,\adState}^\T\mE_i\mR_{i,\adState}
        =\mR_{i,\state}^\T\mE_i\mS_{i,\state}
        =(\mS_{i,\state}^\T\mE_i^\T\mR_{i,\state})^\T.
    \end{equation*}
    Thus, choosing the singular vectors and retained subspaces consistently in \eqref{eqn:SVD} and \eqref{eqn:PBR:proj}, we get
    \begin{align}\label{eqn:comm:proj:mat}
        \begin{aligned}
            \mV_{i,\state} &= \mS_{i,\state}\fU_{i,\state,1}\mathbf{\Sigma}^{-1/2}_{i,\state, 1} = \mR_{i,\adState}\mPsi_{i,\adState,1}\mathbf{\Sigma}^{-1/2}_{i,\adState, 1} = \mW_{i,\adState},\\
            \mW_{i,\state} &= \mR_{i,\state}\mPsi_{i,\state,1}\mathbf{\Sigma}^{-1/2}_{i,\state,1} = \mS_{i,\adState}\fU_{i,\adState,1}\mathbf{\Sigma}^{-1/2}_{i,\adState,1} = \mV_{i,\adState},
        \end{aligned}
    \end{align}
    for $i=1,\ldots,L$.
    Now, let us state system $\switchedSysRed^{\T}$, i.e., the adjoint system of the \ROM \eqref{eq:SLS:ROM}, which reads as
    \begin{equation*}
        \switchedSysRed^{\T} \quad \left\{\quad \begin{aligned}
            -\reduce{\mE}_{\switch(t)}^\T\dot{\adStateRed}(t)  &= \reduce{\mA}^{\T}_{\switch(t)} \adStateRed(t) +  \reduce{\mC}_{\switch(t)}^\T (\outRed(t)-\outVarDes(t)+\tilde{\multiplier}(t)), & \adStateRed(\finalTime) &= \zeroVec, \\
            \reduce{\mE}_{\switch(t^-)}^\T\adStateRed(t^{-})& =\reduce{\mK}^{\T}_{\switch(t^+),\switch(t^-)}\reduce{\mE}_{\switch(t^+)}^\T\adStateRed(t^+),\\
            \adOutRed(t) &= \reduce{\mB}_{\switch(t)}^\T \adStateRed(t),\\
        \end{aligned}\right.
    \end{equation*}
    with $\reduce{\mE}_{i}$, $\reduce{\mA}_{i}$, $\reduce{\mK}_{i,j}$, $\reduce{\mB}_{i}$, and $\reduce{\mC}_{i}$ defined as in \eqref{eqn:SLS:ROM:matrices} via the projection matrices of \eqref{eqn:PBR:proj}.
    The \ROM obtained via \PBR for the adjoint system \eqref{eq:SLS:FOM_Adj} reads as 
    \begin{equation}
        \label{eq:SLS:ROM_Adj}
        \widetilde{\switchedSysAdj} \quad \left\{\quad \begin{aligned}
            -\reduceAd{\mE}_{\switch(t)} \dot{\adStateRedAd}(t)  &= \reduceAd{\mA}_{\switch(t)} \adStateRedAd(t) + \reduceAd{\mC}_{\switch(t)} (\outRed(t)-\outVarDes(t)+\tilde{\multiplier}(t)), & \adStateRedAd(T) &= \zeroVec, \\
            \adStateRedAd(t^-) &= \reduceAd{\mK}_{\switch(t^-),\switch(t^+)} \adStateRedAd(t^+),\\
            \adOutRedAd(t) &= \reduceAd{\mB}_{\switch(t)} \adStateRedAd(t),
        \end{aligned}\right.
    \end{equation}
    where for $i,j=1,\ldots,L$ we have
    \begin{align*}
    	\reduceAd{\mE}_{j} &= \mW_{j,\adState}^\T \mE^{\T}_j \mV_{j,\adState}, &
		\reduceAd{\mC}_{j} &= \mW^{\T}_{j,\adState}\mC_j^{\T}, & 
		\reduceAd{\mA}_{j} &= \mW_{j,\adState}^\T \mA^{\T}_j \mV_{j,\adState}, &
		\reduceAd{\mB}_{j} &= \mB^{\T}_j\mV_{j,\adState},
    \end{align*}
    and $\reduceAd{\mK}_{j,i} = \reduceAd{\mE}_{j}^{-1}\mW_{j,\adState}^{\T}\mE_j^\T\mK^{\mathrm{ad}}_{j,i}\mV_{i,\adState} = \reduceAd{\mE}_{j}^{-1}\mW_{j,\adState}^{\T}\mK_{i,j}^{\T}\mE_i^\T\mV_{i,\adState}$.
    Finally, using \eqref{eqn:comm:proj:mat} we get
    \begin{align*}
        \reduce{\mE}^{\T}_{j} &= \mV^{\T}_{j,\state}\mE^{\T}_j\mW_{j,\state} = \mW_{j,\adState}^{\T}\mE^{\T}_j\mV_{j,\adState} = \reduceAd{\mE}_{j}, &
        \reduce{\mB}^{\T}_{j} &= \mB_{j}^{\T}\mW_{j,\state}\;=\;\mB_j^{\T}\mV_{j,\adState} = \reduceAd{\mB}_{j},\\
        \reduce{\mA}^{\T}_{j} &= \mV^{\T}_{j,\state}\mA^{\T}_j\mW_{j,\state} = \mW_{j,\adState}^{\T}\mA^{\T}_j\mV_{j,\adState} = \reduceAd{\mA}_{j}, &      
        \reduce{\mC}^{\T}_{j} &= \mV_{j,\state}^{\T}\mC_{j}^{\T} = \mW_{j,\adState}^{\T}\mC_j^{\T} = \reduceAd{\mC}_{j},
	\end{align*}
	and
	\begin{align*}
		\reduce{\mE}_j^{-\T}\reduce{\mK}_{i,j}^{\T}\reduce{\mE}_i^\T &= \reduce{\mE}_j^{-\T}
        \left(\reduce{\mE}_i^{-1}\mW_{i,\state}^{\T}\mE_i\mK_{i,j}\mV_{j,\state}\right)^{\T}\reduce{\mE}_i^\T
        = \reduceAd{\mE}_{j}^{-1}\mW_{j,\adState}^{\T}\mK_{i,j}^{\T}\mE_i^\T\mV_{i,\adState}
        = \reduceAd{\mK}_{j,i},
    \end{align*}
    where we used $\mV_{j,\state}=\mW_{j,\adState}$, $\mW_{i,\state}=\mV_{i,\adState}$, and $\reduce{\mE}_j^\T=\reduceAd{\mE}_j$ in the second equality. This concludes the proof of the second statement.
\end{proof}

\begin{remark}
    For \LTI systems in control form, \Cref{teo:duality} holds by replacing the \PBR method with standard \BT.
    For duality in the \FOM, see, e.g., \cite[Sec.~6.2.5]{hinrichsen2026mathematical}.
\end{remark}
\subsection{Application to the \OCP error bounds} \label{subsec:BTbounds}
%
A significant implication of \Cref{teo:duality} is that the projection matrices obtained by applying \PBR to the forward \SLS also reduce its adjoint system. Thus, no separate balancing computation is required for the adjoint, which is why reducing the \OCP is as costly as reducing the \SLS. Particularly, we can use the \PBR error bound for the output of the \FOM and its adjoint to bound their respective outputs in \Cref{lemma:apost_control,lemma:apost_valuefun}.
\begin{corollary}[\PBR error bounds]\label{cor:PBRcontrol}
    In the setting of \Cref{lemma:apost_control,lemma:apost_valuefun}, suppose that
    \eqref{eqn:PBR:bound} applies to the forward and time-reversed adjoint systems.
    Define the primal and dual \PBR estimators
    \begin{equation}\label{eqn:PBR:primal:dual:estimates}
        \begin{aligned}
            \Delta_{\state}\vcentcolon=
            \Delta(\mathbf{\Sigma}_{\state,2},\stateRedOpt,\lambda_1,\inpRedOpt),\qquad
            \Delta_{\adState}\vcentcolon=
            \Delta\bigl(\mathbf{\Sigma}_{\adState,2},\adStateRedOpt,
            \lambda_{1,\adState},\outRedOpt-\outVarDes+\multiplierRedOpt\bigr).
        \end{aligned}
    \end{equation}
    Then,
    \begin{equation}\label{eqn:PBR:value:interval}
        \tilde\eta_-\leq \eta_- \leq
        \calJ_\reg(\outOpt,\inpOpt)-\calJ_\reg(\outRedOpt,\inpRedOpt)
        \leq\eta_+\leq \tilde\eta_+,
    \end{equation}
    where $\eta_-$ and $\eta_+$ are defined in \eqref{eq:value:eta} and
    \begin{equation}\label{eqn:PBR:eta}
        \tilde\eta_\pm\vcentcolon=
        \pm\|\outRedOpt-\outVarDes+\multiplierRedOpt\|_{L^2}\Delta_{\state}
        +\begin{cases}
            -\tfrac{1}{2\lambda}\Delta_{\adState}^2,&\text{for }\tilde\eta_-,\\
            \tfrac{1+\reg}{2}\Delta_{\state}^2,&\text{for }\tilde\eta_+.
        \end{cases}
    \end{equation}
    Moreover, the errors defined in \Cref{lemma:apost_control,lemma:apost_valuefun} satisfy
    \begin{equation}\label{eqn:new:estimate}
        \begin{aligned}
            e_{\outOpt,\inpOpt}&\leq\Delta_{\outOpt,\inpOpt}\leq
            \tilde\Delta_{\outOpt,\inpOpt}(\outRedOpt,\inpRedOpt,\multiplierRedOpt)
            \vcentcolon=
            \left(\tfrac1\lambda\Delta_{\adState}^2
            +(1+\reg)\Delta_{\state}^2\right)^{1/2},\\
            e_{\calJ_\reg}&\leq\Delta_{\calJ_\reg}\leq
            \tilde\Delta_{\calJ_\reg}(\outRedOpt,\inpRedOpt,\multiplierRedOpt)
            \vcentcolon=\max\{|\tilde\eta_-|,|\tilde\eta_+|\},
        \end{aligned}
    \end{equation}
    where the estimators $\Delta_{\outOpt,\inpOpt}$ and
    $\Delta_{\calJ_\reg}$ are defined in
    \eqref{eq:abstract_control_bound} and \eqref{eq:value:convex:bound}, respectively.
\end{corollary}
\begin{proof}
    Applying \eqref{eqn:PBR:bound} to the time-reversed adjoint system gives,
    for any $\out,\multiplier\in L^2(0,T;\R^\outDim)$,
    \begin{equation}\label{eqn:PBR:adjoint:bound}
        \|\adOut(\out,\multiplier)-\adOutRed(\out,\multiplier)\|_{L^2}
        \leq
        \Delta\bigl(\mathbf{\Sigma}_{\adState,2},
        \adStateRed(\out,\multiplier),\lambda_{1,\adState},
        \out-\outVarDes+\multiplier\bigr),
    \end{equation}
    By \Cref{teo:duality} and
    \eqref{eq:OptCond:ROMMY:Adstate1}--\eqref{eq:OptCond:ROMMY:Adstate3},
    the reduced adjoint driven by $\outRedOpt-\outVarDes+\multiplierRedOpt$
    has state $\adStateRedOpt$ and output $\adOutRedOpt$.
    Thus,
    \begin{equation*}
        \|\out(\inpRedOpt)-\outRedOpt\|_{L^2}\leq\Delta_{\state},
        \qquad
        \|\adOut(\outRedOpt,\multiplierRedOpt)-\adOutRedOpt\|_{L^2}
        \leq\Delta_{\adState}.
    \end{equation*}
    Inserting these inequalities into \eqref{eq:abstract_control_bound}
    proves the control and output estimate by \Cref{lemma:apost_control}.
    For the lower bound in \eqref{eqn:PBR:eta},
    the Cauchy--Schwarz inequality and the primal \PBR estimate give
    \begin{equation*}
        \bigl\langle\outRedOpt-\outVarDes+\multiplierRedOpt,
        \out(\inpRedOpt)-\outRedOpt\bigr\rangle_{L^2}
        \geq-\|\outRedOpt-\outVarDes+\multiplierRedOpt\|_{L^2}\Delta_{\state}.
    \end{equation*}
    Together with the dual \PBR estimate, this gives
    $\tilde\eta_-\leq\eta_-$ by \eqref{eq:value:eta},
    and \Cref{lemma:apost_valuefun} proves the lower bound
    in \eqref{eqn:PBR:value:interval}.
    Turning to the upper bound, Cauchy--Schwarz bounds the same inner product from above
    by $\|\outRedOpt-\outVarDes+\multiplierRedOpt\|_{L^2}\Delta_{\state}$.
    Together with the primal \PBR estimate, this gives
    $\eta_+\leq\tilde\eta_+$ and proves the upper bound.
\end{proof}

\begin{remark}
	\Cref{cor:PBRcontrol} is based on the a posteriori error bound \eqref{eqn:PBR:bound}, which is associated with the specific \MOR procedure described in \Cref{subsec:PBR}. We emphasize that, more generally, any reduction methodology that is commutative with the optimality conditions for \SLS, i.e., for which \Cref{teo:duality}\,\ref{item2} holds and that is equipped with a certified error estimator can be used to establish corollaries of this type. Clearly, the computational cost of evaluating such an error bound should be favorable compared with the cost of evaluating the \FOM at the optimal input provided by the~\ROM-\OCP.
\end{remark}
\subsection{Computational procedure and uniform operator approximation} \label{subsec:computational_operator_approx}
%
The offline--online computational procedure can be summarized as follows.

\subsubsection{Offline phase}
Solve the \GLEs in \eqref{eqn:GLE:b}, construct the balancing bases via \eqref{eqn:PBR:proj}, build a hierarchy of \ROMs for selected dimensions $r\leq n$, and prepare the data for the error estimator~\eqref{eqn:PBR:bound}.
By \Cref{teo:duality}, the forward bases also determine the reduced
adjoint, so no separate balancing computation is required.
The discarded singular values coincide for both systems, but the
adjoint trajectory and \LMI-defect contributions must be treated separately.
In the special case of \eqref{eqn:PBR:special:bound}, precompute $c_r$. For fixed system matrices, the projection bases and, in this special case, $c_r$ are independent of the parameter tuple $\bmu$.

\subsubsection{Online phase}
For the current parameter tuple $\bmu$, solve the \ROM-\OCP and evaluate
the estimators in \Cref{cor:PBRcontrol}.
If the certificate exceeds the prescribed tolerance, increase $r$
and repeat the solve with a larger precomputed \ROM.
In the special case of \eqref{eqn:PBR:special:bound}, we have the online-efficient choice
\begin{equation}\label{eq:est_cheap}
    \tilde\Delta_{\outOpt,\inpOpt}(\outRedOpt,\inpRedOpt,\multiplierRedOpt)
    =c_r\left(\tfrac{1}{\lambda}
    \|\outRedOpt-\outVarDes+\multiplierRedOpt\|_{L^2}^2
    +(1+\reg)\|\inpRedOpt\|_{L^2}^2\right)^{1/2}.
\end{equation}
With this bound, or a certified offline--online decomposition in the
general \PBR case, the procedure is online-adaptive and online-efficient.
These are online-efficient \PBR-based bounds, in contrast to
the full-order-evaluation bounds in
\Cref{lemma:apost_control,lemma:apost_valuefun}.
In this case, the neglected singular-value tail controls the uniform approximation of \eqref{eq:sol_op}.
\begin{corollary}[Uniform approximation of the solution operator]\label{rem:uniform:solution}
    Let $\bar{U}\subseteq \R^\inpDim$ and $\bar{Y}\subseteq \R^\outDim$ be nonempty and compact,  $0<\lambda_{\min}\leq \lambda_{\max}$,  $0<\gamma_{\min}\leq \gamma_{\max}$, and $\bar{\finalTime} >0$.
    In addition to the standing assumptions on $\Uad$ and $\Yad$, define $\calM$ as the set of parameter tuples $\bmu$ from \eqref{eq:definition_mu} satisfying
    $\switch\in\calP$, $0<\finalTime\leq\bar{\finalTime}$,
    $\lambda\in[\lambda_{\min},\lambda_{\max}]$, $\reg\in[\reg_{\min},\reg_{\max}]$,
    $\outVarDes\in L^2(0,\finalTime;\bar{Y})$, $\inpVarDes\in L^2(0,\finalTime;\bar{U})$,
    and $\Uad\subseteq\bar{U}$, $\Yad\subseteq\bar{Y}$.
    Further, assume that \eqref{eqn:PBR:special:bound} holds uniformly, with
    $c_r=2\operatorname{tr}(\mathbf{\Sigma}_{1,\state,2})$ independent of
    $\bmu\in\mathcal M$.
    Then,
    \begin{equation*}
        \sup_{\bmu\in\mathcal M}
        \tilde\Delta_{\outOpt,\inpOpt}(\bmu)
        \leq C c_r,
    \end{equation*}
    where $C$ is independent of $\bmu$ and $r$.
    Hence, with $\Solper_r(\bmu)
    \vcentcolon=(\outRedOpt(\bmu),\inpRedOpt(\bmu))$,
    \begin{equation*}
        \sup_{\bmu\in\mathcal M}
        \|\Solper(\bmu)-\Solper_r(\bmu)\|_{L^2\times L^2}
        \leq Cc_r \to 0 \qquad \text{for $r\to n$}.
    \end{equation*}

\end{corollary}
\begin{proof}
    The uniformly bounded time horizons, admissible
    input sets and output references, and uniformly bounded reduced
    input-to-output operators imply
    \begin{equation*}
        \|\inpRedOpt(\bmu)\|_{L^2}\leq M_\inp,\qquad
        \|\outRedOpt(\bmu)-\outVarDes+\multiplierRedOpt(\bmu)\|_{L^2}\leq M_\out.
    \end{equation*}
    Inserting these bounds into \eqref{eq:est_cheap} gives
    \begin{equation*}
        \tilde\Delta_{\outOpt,\inpOpt}(\bmu)\leq C_0c_r,
        \qquad
        C_0^2=\tfrac{M_\out^2}{\lambda_{\min}}+(1+\reg_{\max})M_\inp^2.
    \end{equation*}
    Using \eqref{eq:err:OCP:ROM},
    $\lambda\geq\lambda_{\min}$, and \Cref{cor:PBRcontrol} it follows
    \begin{equation*}
        \begin{aligned}
            \|\Solper(\bmu)-\Solper_r(\bmu)\|_{L^2\times L^2}^2
            &=\|\outOpt-\outRedOpt\|_{L^2}^2+\|\inpOpt-\inpRedOpt\|_{L^2}^2
            \leq \max\{1,\lambda_{\min}^{-1}\}\,C_0^2c_r^2.
        \end{aligned}
    \end{equation*}
    Thus, both claims hold with
    $C=\max\{1,\lambda_{\min}^{-1/2}\}\,C_0$,
    which is independent of $\bmu$ and $r$.
\end{proof}

\section{Numerical Experiments}\label{sec:Numerics}
%
The purpose of this section is to provide a numerical validation of the a posteriori bounds established in \Cref{sec:outputbased_error_est,sec:BT_for_OCP}, and to evaluate both the approximation accuracy and the computational speed-ups achieved.
All simulations were performed with MATLAB R2024b on a 2020 MacBook Pro with a 2.3 GHz quad-core Intel Core i7 processor and 16 GB of RAM. The source code is accessible here\footnote{Data will be made publicly available upon publication.}.

\subsection{Computational details}
%
We choose the time interval $[0,\finalTime]=[0,1]$ and discretize it using the implicit Euler method with $K=500$ time steps, whereas the cost functional is discretized by a right-hand rectangle rule, making time discretization and optimization commute.
As a solution method for the \OCP, we use \MY continuation with increased penalties $\gamma_k=\min(2^{k-1},\gamma_{\max})$. Increasing $\gamma_k$ progressively strengthens the penalization of violations of the output constraints. As an inner solver for fixed $\gamma_k$, we use a Barzilai--Borwein proximal gradient method \cite{azmi2023nonmonotone}. The estimators are only evaluated once, after the last inner \ROM \OCP solve for the last $\gamma_k$. The parameters and stopping tolerances are specified below. We impose the pointwise box constraints $\inp_a(t)\leq\inp(t)\leq\inp_b(t)$ and $\out_a(t)\leq\out(t)\leq\out_b(t)$.
The reported online speed-up is the \FOM solution time divided by the sum of the \ROM solution and estimator evaluation times, averaged over 5 repetitions.

%
\subsection{Synthetic example}\label{subsec:Numerics:synthetic}

Consider a two-mode \SLS of dimension~$\stateDim=10^4$ with state matrices $\mA_1,\mA_2\in\R^{\stateDim\times\stateDim}$ defined by
\begin{align}\label{eqn:synt:ex}
    \begin{aligned}
        \mA_1(i,j)\vcentcolon=\begin{cases} -1, & \mbox{if }i=j, \\
            \tfrac{1}{2}, & \mbox{if }i-1=j,\\
            0,&\mbox{otherwise},
        \end{cases}\quad\quad \mA_2(i,j)\vcentcolon=\begin{cases} -2, & \mbox{if }i=j, \\
            \tfrac{4}{5}, & \mbox{if }i-1=j,\\
            -\tfrac{1}{5}, & \mbox{if }i+1=j,\\
            0,&\mbox{otherwise}.
        \end{cases}
    \end{aligned}
\end{align}
The remaining matrices are $\mE_1=\mE_2=\mI_\stateDim$, $\mB_1^\T=\mC_1=[0,\ldots,0,1]$, $\mB_2^\T=\mC_2=[1,0,\ldots,0]$, and $\mK_{i,j}=\mathbf{I}_{\stateDim\times\stateDim}$. The system matrices are taken from \cite{PonGB20}. We set $\inpVarDes\equiv0$ and to construct the \BT-\ROM, we use the stopping tolerance $\varepsilon=10^{-11}$ to compute the solution of \GLEs, and precompute the \ROM offline for $r= 2,\ldots,16$. Online, we treat the following two set-ups of $\bmu$ introduced in \eqref{eq:definition_mu}:
\begin{itemize}
    \item Configuration $\bmu_1$: we choose the target $\outVarDes=\out(\tilde\inp)$ with $\tilde\inp\equiv1$ and set $\lambda=10^{-1}$.
    The bounds are $\inp_a\equiv0.1$, $\inp_b\equiv0.6$, $\out_a\equiv0$, and $\out_b\equiv0.1$. We use $N = 5$ randomly chosen switches starting in mode~$1$, and $10^{-14}$ as optimization tolerance.
    \item Configuration $\bmu_2$: we choose the target $\outVarDes(t)=2+\sin(2\pi t)$ and set $\lambda=10^{-3}$.
    The bounds are $\inp_a\equiv0.01$, $\inp_b\equiv20$, $\out_a\equiv0$, and $\out_b\equiv2$. We use $N = 20$ randomly chosen switches starting in mode~$2$. The inner optimization tolerance is chosen as $10^{-13}$.
\end{itemize}
%

\begin{figure}[t]
    \centering
    \begin{subfigure}[t]{0.48\textwidth}
        \centering
        \includegraphics{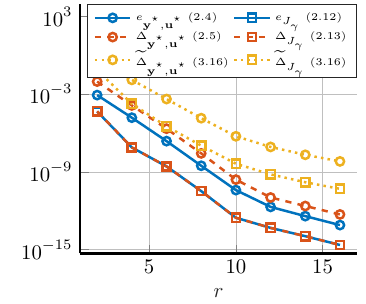}
        \subcaption{Errors and associated estimators.}
        \label{fig:toy_est}
    \end{subfigure}
    \hfill
    \begin{subfigure}[t]{0.48\textwidth}
        \centering
        \includegraphics{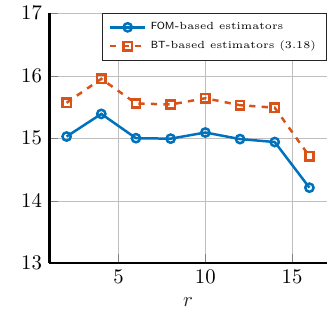}
        \subcaption{Averaged speed-up including solution of the \ROM-\OCP and evaluation of the error estimators.}
        \label{fig:toy_speedup}
    \end{subfigure}
    \begin{subfigure}[t]{0.48\textwidth}
        \centering
        \includegraphics{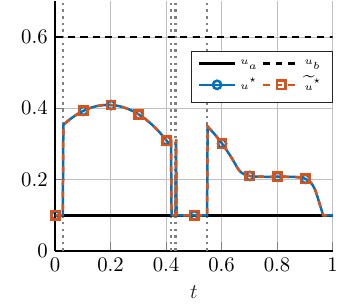}
        \subcaption{Optimal control and associated constraints. Switching times are indicated by dotted lines.}
        \label{fig:toy_inp:con}
    \end{subfigure}
    \hfill
    \begin{subfigure}[t]{0.48\textwidth}
        \centering
        \includegraphics{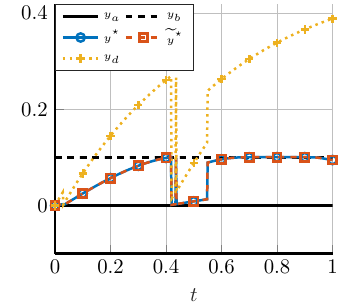}
        \subcaption{Optimal output and associated constraints.}
        \label{fig:toy_out:con}
    \end{subfigure}
    \caption{Synthetic example, configuration $\bmu_1$ with \MY regularization $\reg=250$. The lower panels are for $\stateDimRed=8$. Gray dotted lines in the control panel mark the switching times.}
    \label{fig:synthetic_ex_config_1}
\end{figure}

\begin{figure}[t]
    \centering
    \begin{subfigure}[t]{0.48\textwidth}
        \centering
        \includegraphics{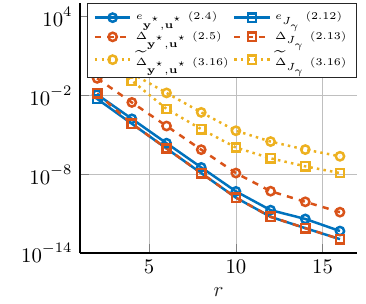}
        \subcaption{Errors and associated estimators.}
\label{fig:toy_est:2}
    \end{subfigure}
    \hfill
    \begin{subfigure}[t]{0.48\textwidth}
        \centering
        \includegraphics{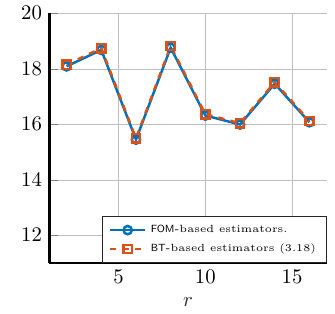}
        \subcaption{Averaged speed-up including solution of the \ROM-\OCP and evaluation of the error estimators.}
\label{fig:toy_speedup:2}
    \end{subfigure}
    \\
    \begin{subfigure}[t]{0.48\textwidth}
        \centering
        \includegraphics{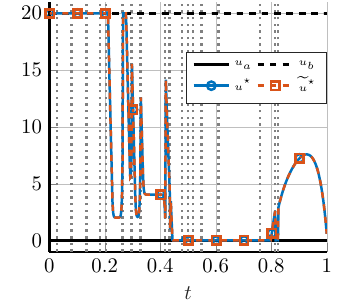}
        \subcaption{Optimal control and associated constraints. Switching times are indicated by dotted lines.}
        \label{fig:toy_inp:con:2}
    \end{subfigure}
    \hfill
    \begin{subfigure}[t]{0.48\textwidth}
        \centering
        \includegraphics{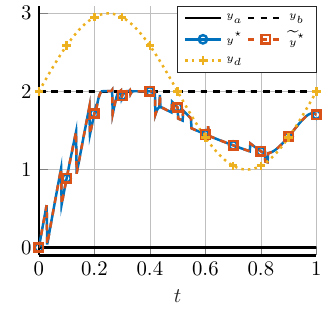}
        \subcaption{Optimal output and associated constraints. Note that the jumps in the output are precisely at the switching times.}
        \label{fig:toy_out:con:2}

    \end{subfigure}
    \caption{Synthetic example, configuration $\bmu_2$ with \MY regularization $\reg=250$. The lower panels are for $\stateDimRed=8$. Gray dotted lines in the control panel mark the switching times.}
    \label{fig:synthetic_ex_config_2}
\end{figure}

The results for $\bmu_1$ and $\bmu_2$ are shown in \Cref{fig:synthetic_ex_config_1,fig:synthetic_ex_config_2}, respectively.
As $\stateDimRed$ increases from $2$ to $16$, the combined error decreases from $8.6\cdot10^{-4}$ to $7.8\cdot10^{-14}$ for $\bmu_1$, and from $1.2\cdot10^{-2}$ to $4.9\cdot10^{-13}$ for $\bmu_2$. 
\Cref{fig:toy_est,fig:toy_est:2} compare the computed errors $e_{\calJ_\reg}$ and $e_{\outOpt,\inpOpt}$ with the full-order evaluation estimators from \Cref{lemma:apost_control,lemma:apost_valuefun} and the \BT-based indicator obtained from the expressions in \Cref{cor:PBRcontrol}, using the simplified error estimate \eqref{eqn:PBR:special:bound}, as this is applicable for this problem, as shown in \cite[Sec.~5.1]{ManU26}. The \FOM output-constraint violation $\|\outOpt-\Proj(\outOpt)\|_{L^2}$ is approximately $3.80\cdot10^{-4}$ for $\bmu_1$ and $1.85\cdot10^{-3}$ for $\bmu_2$, both at $\reg=250$, which is negligible and, correspondingly, not discernible in \Cref{fig:toy_out:con,fig:toy_out:con:2}.

\Cref{fig:toy_est,fig:toy_est:2} also show that the full-order-evaluation estimator for the combined error exceeds the computed error by factors of $5.3$--$11.1$ for $\bmu_1$ and $16.9$--$27.6$ for $\bmu_2$.
The value-function estimator almost coincides with the computed cost error, except at $r=2$ for $\bmu_2$, where the ratio is $2.63$. At the smallest cost errors, some computed estimates fall below the measured errors by less than $10^{-16}$ in absolute terms, indicating the limit of numerical resolution.
The \BT-based estimators are larger than the corresponding errors by factors ranging from approximately $4.5\cdot10^2$ to $4.9\cdot10^5$ across both configurations. However, the ratios between the mean time needed for the evaluation of the \FOM-based estimator and the mean time needed for the evaluation of the \BT-based estimator, averaged over all dimensions and timing repetitions, are approximately $3000$ for both configurations. Thus, while the \BT-based error estimate is less sharp than its \FOM-based counterpart, it is substantially faster to evaluate. It should be noted that, when the solution time of the \ROM-\OCP is included, the resulting computational advantage is substantially reduced. This is because the time required to solve the \ROM-\OCP dominates the online computational cost and, consequently, also exceeds the time needed to evaluate both the \FOM-based and \BT-based error estimators, particularly for $\bmu_2$; see \Cref{fig:toy_speedup,fig:toy_speedup:2}. Increasing the \FOM size would, as expected, raise the computational burden of the \FOM-based estimator, whereas it would essentially not affect the computational requirements of the \BT-based estimator. The nonmonotone speed-up curve for $\bmu_2$ reflects the varying numbers of Barzilai--Borwein iterations across reduced dimensions.
%
%

\subsection{Sparse economic optimal control of the advection-diffusion-reaction equation}
%
We consider the two-room heating example from \cite[Sec.~2.3 and 2.5]{KarMUV26}, which is inspired by \cite{SchU18}.
The rectangular domain $\Omega\subseteq\R^2$ consists of the rooms $\Omega_1$ and $\Omega_3$, the door $\Omega_2$, and the walls $\Omega_4$ (see \cite[Fig.~1]{KarMUV26} for details). Its boundary is partitioned into two disjoint sets $\Gamma_c$ and $\Gamma_o$.
The control $\inp\colon[0,1)\to\R^{\inpDim}$ acts on $\Gamma_c$ through $\inpDim=10$ actuators.
The switching signal $\switch\colon[0,1)\to\{1,2\}$ indicates whether the door is closed ($\switch=1$) or open ($\switch=2$).
We fix two switches at $t=0.028$ and $t=0.436$, with mode sequence $2\to1\to2$. The goal is to control the following switched heat equation
\begin{subequations}
    \label{eq: parabolic_problem}
    \begin{align}
        \label{eq: parabolic_problem-a}
        \zeta(\bx,\switch(t))\partial_t \StateVar(t,\bx) + \calA(\bx,\switch(t)) \StateVar(t,\bx)&= 0, && (t,\bx) \in (0,1)\times\Omega,\\
        \kappa(\bs,\switch(t)) \tfrac{\partial}{\partial n}\StateVar(t,\bs) + \gamma_{\mathrm{o}}\StateVar(t,\bs)&= 0, && (t,\bs) \in (0,1) \times \Gamma_{\mathrm{o}},\\
        \kappa(\bs,\switch(t)) \tfrac{\partial}{\partial n}\StateVar(t,\bs) - \sum\nolimits_{i=1}^{\inpDim} \chi_i(\bs) \inp_i(t)&=0, && (t,\bs) \in (0,1) \times \Gamma_{\mathrm{c}},\\
        \theta(0,\bx) &= \theta_0(\bx), && \bx \in \Omega,
    \end{align}
\end{subequations}
with the differential operator
\begin{align*}
    \calA(\bx,\switch)\StateVar(t,\bx)\vcentcolon=-\nabla \cdot\big(\kappa(\bx,\sigma)\nabla \StateVar(t,\bx)\big) + b(\bx)\cdot\nabla\StateVar(t,\bx) + c(\bx)\StateVar( t,\bx).
\end{align*}
The coefficients $\zeta$ and $\kappa$ are piecewise constant, and switching occurs only on the door $\Omega_2$ (see \cite[Sec.~4.2]{KarMUV26} for the precise definitions and parameter values).
%
%
The output consists of the average temperatures in the two rooms $\Omega_1$ and $\Omega_3$.
Accordingly, we set $\outDim=2$ and define $\calC:L^2(\Omega)\to\R^{\outDim}$ by $\varphi\mapsto y=(y_1,y_2)=\calC\varphi$, where $y_i\vcentcolon=|D_i|^{-1}\int_{D_i}\varphi(\bx)\,\mathrm d\bx$ for $D_1=\Omega_1$ and $D_2=\Omega_3$.
To promote sparse controls, we add an $\ell^1$-term to the objective.
For $c_0\geq0$, $c_1>0$, and a fixed shift $\shift$ (cf.~\Cref{rem:AL:shifts}), the computed objective is
\begin{align*}
    \mathcal J_{\reg,\shift}(\out,\inp) = \tfrac{c_0}{2}\|\out-\outVarDes\|^2_{L^2}
  	+& \tfrac{\lambda}{2}\|\inp\|^2_{L^2}
    + c_1\int_0^1\|\inp(t)\|_{\ell^1}\,\dt\\
    +& \tfrac{\reg}{2}\|\out+ \shift/\reg-\Proj(\out+ \shift/\reg)\|_{L^2}^2.
\end{align*}
The argument for convex control regularization in the proof of \cite[Thm.~3]{KarMUV26} extends the error analysis from \Cref{lemma:apost_control} to the additional $\ell^1$-term. The value-function interval in \Cref{lemma:apost_valuefun} is enlarged by the Lipschitz constant $c_1\sqrt{\inpDim}$ of the $\ell^1$-norm term multiplied by the control-error bound, since $T=1$.

Piecewise linear finite elements yield a discrete system of dimension $\stateDim=2170$ with nonsingular mass matrices $\mE_j$ and identity transition matrices $\mK_{i,j}$.
We compute the \PBR models from \Cref{sec:BT_for_OCP} as follows:
The generalized Lyapunov equations are solved by a low-rank iteration with LR-ADI inner solves with a tolerance of $10^{-10}$ for the residual norm. The \PBR-\ROMs for each $r=10,\ldots,60$ are precomputed offline, and tested online for the following configurations:
\begin{itemize}
    \item Configuration $\bmu_3$: a staggered room-heating schedule with $\theta_0=0$.
    Room 1 is brought to its operating temperature first, and room 2 is heated for later use after the door reopens.
    For $S_{a,b}(t)=3z^2-2z^3$, where $z=\min\{1,\max\{0,(t-a)/(b-a)\}\}$, we prescribe
    $y_{d,1}(t)=\tfrac{9}{10}S_{0,0.15}(t)$ and $y_{d,2}(t)=S_{0.45,0.65}(t)$.
    We use $c_0=100$, $\lambda=10^{-1}$, and $c_1=10^{-4}$. The boundary heaters satisfy $0\leq u_i\leq3$, and the output bounds are $\out_a=(0,0)^\T$ and $\out_b=(0.95,1.1)^\T$.
    \item Configuration $\bmu_4$: economic heating with staggered occupancy, $\theta_0\equiv3$, $c_0=0$, $\lambda=10^{-1}$, and $c_1=10^{-4}$.
    We minimize heating effort subject to $0\leq u_i\leq3$, $\out_b\equiv(3,3)^\T$, and the time-dependent lower bounds
    $y_{a,1}(t)=\tfrac{29}{10}-\tfrac{4}{10}S_{0.60,0.80}(t)$ and $y_{a,2}(t)=\tfrac{23}{10}+\tfrac{4}{10} S_{0.50,0.75}(t)$.
    The initial-state solution is subtracted from the bounds to obtain a zero-initial-state optimization problem and added back for plotting.
\end{itemize}
For each configuration, an augmented-Lagrangian iteration (see \Cref{rem:AL:shifts}) on a fixed $\stateDimRed=20$ \ROM determines $\shift$, starting from the zero control and multiplier and using $\gamma_k=\min(2^{k-1},1000)$. We retain the shift used in its final subproblem. Every \FOM and \ROM comparison then solves one shifted \MY problem with this same fixed $\shift$ and $\gamma=1000$, as in \Cref{rem:AL:shifts}.

For each optimal \ROM control, the \PBR approximation solves one enriched primal and adjoint \ROM of dimension $\stateDimRed=81$ and computes their output differences to the optimization \ROM. Insertion into the \OCP estimates, including the $\ell^1$-norm cost correction, gives $\Delta^{\mathrm{\PBR}}_{\outOpt,\inpOpt}$ and $\Delta^{\mathrm{\PBR}}_{J_\gamma}$. The online calculation uses no \FOM solve and has a state dimension independent of $\stateDim$. Again, we remark that approximating the trajectories and truncating the hierarchy at $\tilde \stateDimRed=90<\stateDim$ makes this an error indicator rather than a certified bound. The optimization tolerance is $10^{-10}$.

\begin{figure}[t]
    \centering
    \begin{subfigure}[t]{0.48\textwidth}
        \centering
        \includegraphics{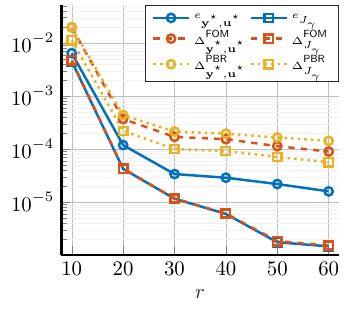}
        \subcaption{Errors and associated estimators.}
    \end{subfigure}
    \hfill
    \begin{subfigure}[t]{0.48\textwidth}
        \centering
        \includegraphics{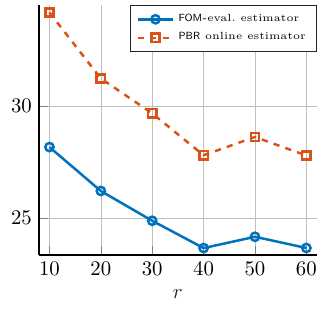}
        \subcaption{Averaged speed-up including solution of the \ROM-\OCP and evaluation of the error estimators.}
    \end{subfigure}
    \\
    \begin{subfigure}[t]{0.48\textwidth}
        \centering
        \includegraphics{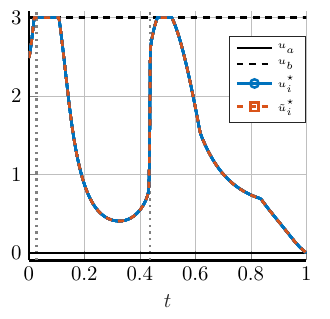}
        \subcaption{Optimal control and associated constraints. Switching times are indicated by dotted lines.}
    \end{subfigure}
    \hfill
    \begin{subfigure}[t]{0.48\textwidth}
        \centering
        \includegraphics{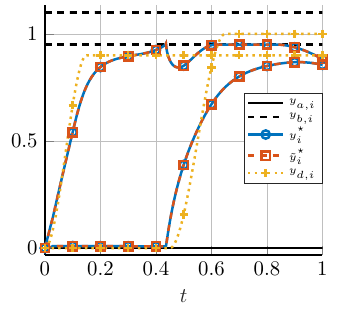}
        \subcaption{Optimal output and associated constraints.}
    \end{subfigure}
    \caption{Two-room heating example with $\stateDim=2170$. Configuration $\bmu_3$: staggered temperature targets, nonnegative bounded heating, and an overheating limit in room 1. The \FOM and \ROM curves solve shifted \MY problems at $\gamma=1000$ with a common fixed shift from an $\stateDimRed=20$ \ROM. Gray dotted lines in the control panel mark the switching times.}
    \label{fig:heat2d_ex_config_3}
\end{figure}

\begin{figure}[t]
    \centering
    \begin{subfigure}[t]{0.48\textwidth}
        \centering
        \includegraphics{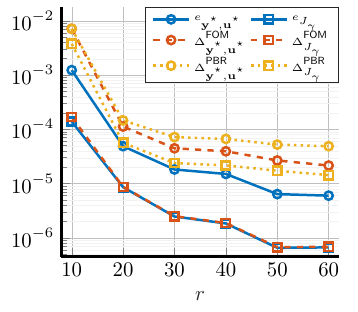}
        \subcaption{Errors and associated estimators.}
    \end{subfigure}
    \hfill
    \begin{subfigure}[t]{0.48\textwidth}
        \centering
        \includegraphics{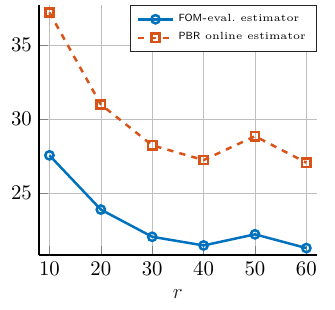}
        \subcaption{Averaged speed-up including solution of the \ROM-\OCP and evaluation of the error estimators.}
    \end{subfigure}
    \\
    \begin{subfigure}[t]{0.48\textwidth}
        \centering
        \includegraphics{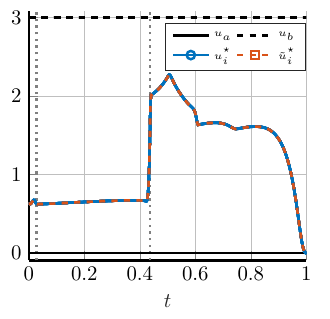}
        \subcaption{Optimal control and associated constraints. Switching times are indicated by dotted lines.}
    \end{subfigure}
    \hfill
    \begin{subfigure}[t]{0.48\textwidth}
        \centering
        \includegraphics{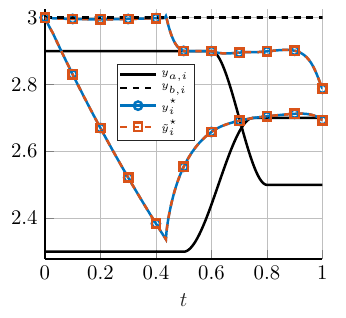}
        \subcaption{Optimal output and associated constraints.}
    \end{subfigure}
    \caption{Two-room heating example with $\stateDim=2170$.  Configuration $\bmu_4$: minimum heating effort with staggered occupancy, time-dependent comfort bounds, and stored initial heat. The \FOM and \ROM curves solve shifted \MY problems at $\gamma=1000$ with a common fixed shift from an $\stateDimRed=20$ \ROM. Gray dotted lines in the control panel mark the switching times.}
    \label{fig:heat2d_ex_config_4}
\end{figure}

\Cref{fig:heat2d_ex_config_3,fig:heat2d_ex_config_4} display errors, estimators, speed-ups, and the $\stateDimRed=60$ trajectories, analogously to the synthetic example. A single \FOM solution per configuration is the reference for all reduced dimensions. For $\bmu_4$ with $c_0=0$, the combined error contains only the control and multiplier errors, with no direct output-error contribution.

For $\bmu_3$, the combined error decreases from $6.43\cdot10^{-3}$ at $\stateDimRed=10$ to $1.60\cdot10^{-5}$ at $\stateDimRed=60$. The \FOM-evaluated estimates are $1.97\cdot10^{-2}$ and $9.01\cdot10^{-5}$, while the \PBR approximations are $1.97\cdot10^{-2}$ and $1.43\cdot10^{-4}$.
For $\bmu_4$, the combined error decreases from $1.23\cdot10^{-3}$ at $r=10$ to $5.97\cdot10^{-6}$ at $r=60$. The \FOM-evaluated estimates range from $7.12\cdot10^{-3}$ to $2.15\cdot10^{-5}$, while the \PBR approximations range from $7.14\cdot10^{-3}$ to $4.87\cdot10^{-5}$.
The \PBR approximation retains a visible tail floor at the larger reduced dimensions. Both estimators exceed the measured combined and cost errors at every tested dimension.

The evaluation of the \PBR error indicator takes about $16.8\,\mathrm{ms}$ for $\bmu_3$ and $15.3\,\mathrm{ms}$ for $\bmu_4$, versus $368.9\,\mathrm{ms}$ and $347.4\,\mathrm{ms}$ for the error estimate based on the \FOM evaluation with the optimal control and optimal output derived from the \PBR-\ROM. Including the \ROM solve, the online speed-ups are $27.8$--$34.2$ (\PBR-based error indicator) versus $23.7$--$28.2$ (\FOM-based error estimate) for $\bmu_3$, and $27.1$--$37.2$ versus $21.3$--$27.6$ for $\bmu_4$.
These online times exclude system projection, \PBR preparation, shift generation, and the common Gramians computation, which takes $60.9\,\mathrm{s}$ with four computation threads.

The control norms are $\|\inpOpt\|_{L^2}=5.3447$ for $\bmu_3$ and $4.1925$ for $\bmu_4$. With activity threshold $|\multiplier_i^\star|>10^{-8}$, the first-output \MY multiplier in $\bmu_3$ is active in $23.2\%$ of time steps, while heaters reach their upper limit in $14.2\%$ of actuator-time samples, and $\|\outOpt-\outVarDes\|_{L^2}=0.1388$. In $\bmu_4$, the lower-bound \MY multipliers are active in $(10.6\%,5.0\%)$ of time steps and the first-room upper-bound multiplier in $3.6\%$, while the heater upper limit remains inactive.
Both hard-constrained discretized \OCPs are feasible: zero control suffices for $\bmu_3$. An independent $\bmu_4$ witness satisfies the input bounds and has output margin at least $2.45\cdot10^{-3}$ at all positive grid times, verified using the original finite-element matrices and the \textsc{Matlab} \FOM; the upper bound is attained at $t=0$. The shifted \MY-\FOM optima have maximum output-constraint violations $5.27\cdot10^{-3}$ and $7.05\cdot10^{-3}$, respectively.

\section{Conclusion}
We developed a balancing-based \ROM for \MY-regularized, output-input-constrained optimal control of \SLS with an unknown switching signal in the set \eqref{eq:suitableSwitchingSignals} and potential state transition jumps at the switching times.
By combining a posteriori error bounds for the optimal control, optimal output, and optimal value function, namely \Cref{lemma:apost_control,lemma:apost_valuefun}, we obtain a certified approximation of the regularized \OCP solution.

Subsequently, we propose employing the \PBR-\MOR framework as a \ROM surrogate for the \SLS. The resulting procedure incorporates an error estimate that can be leveraged to derive computationally efficient error indicators for the \OCP; see \Cref{subsec:PBR} and \Cref{cor:PBRcontrol}. A key ingredient that allows this is a duality (or commutativity) result, which ensures that the reduced forward and adjoint systems can be obtained via the same balancing construction; see \Cref{teo:duality}. 

Overall, by exploiting the offline–online decomposition of the error estimators and indicators, the proposed approach yields an online-efficient procedure for the reliable solution of the \ROM-\OCP.
%
\bibliographystyle{siamplain}
\bibliography{journalabbr.bib,biblio.bib}

@book{HPUU09,
	AUTHOR = {Hinze, M. and Pinnau, R. and Ulbrich, M. and Ulbrich, S.},
	TITLE = {Optimization with {PDE} constraints},
	SERIES = {Mathematical Modelling: Theory and Applications},
	VOLUME = {23},
	PUBLISHER = {Springer},
	ADDRESS = {Dordrecht},
	YEAR = {2009},
	PAGES = {xii+270},
	DOI = {10.1007/978-1-4020-8839-1}
}

@book{Ant09,
	title = {Approximation of Large-Scale Dynamical Systems},
	author={Antoulas, A. C.},
	publisher={SIAM},
	series ={Advances in Design and Control},
	doi = {10.1137/1.9780898718713},
	year ={2005},
	address = {Philadelphia, PA}
}

@book{NMT11,
	title = {Reduced-Order Modelling for Flow Control},
	address = {Vienna},
	series = {CISM International Centre for Mechanical Sciences},
	doi = {10.1007/978-3-7091-0758-4},
	year = {2011},
	publisher = {Springer},
	editor = {Noack, B. R. and Morzy{\'n}ski, M. and Tadmor, G.},
	volume = {528}
}

@article{TV09,
	author={Tr\"oltzsch, F. and Volkwein, S.},
	title={{POD} a-posteriori error estimates for linear-quadratic optimal control problems},
	journal={Comput. Optim. Appl.},
	volume={44},
	pages={83--115},
	doi={10.1007/s10589-008-9224-3},
	year={2009}}

@article{ITO2003221,
title = {Semi-smooth {Newton} methods for state-constrained optimal control problems},
journal = SCL,
volume = {50},
number = {3},
pages = {221--228},
year = {2003},
issn = {0167-6911},
doi = {10.1016/S0167-6911(03)00156-7},
author = {Kazufumi Ito and Karl Kunisch},
}

@article{HJ18,
	author = {Heinkenschloss, M. and Jando, D.},
	title = {Reduced order modeling for time-dependent optimization problems with initial value controls},
	journal = SIAMSciComp,
	volume = {40},
	pages = {A22--A51},
	year = {2018},
	doi = {10.1137/16M1109084},
	number = {1}
}

@incollection{GV17,
	author = {Gubisch, M. and Volkwein, S.},
	title = {Proper Orthogonal Decomposition for Linear-Quadratic Optimal Control},
	booktitle = {Model Reduction and Approximation: Theory and Algorithms},
	pages = {3--63},
	year={2017},
	chapter = {1},
	publisher={SIAM},
	address={Philadelphia},
	doi = {10.1137/1.9781611974829.ch1},
	editor = {Benner, P. and Ohlberger, M. and Cohen, A. and Willcox, K.},
	series = {Computational Science \& Engineering},
	volume = {15}
}

@article{SBMR18,
	author = {Strazzullo, M. and Ballarin, F. and Mosetti, R. and Rozza, G.},
	title = {Model reduction for parametrized optimal control problems in environmental marine sciences and engineering},
	journal = SIAMSciComp,
	volume = {40},
	pages = {B1055--B1079},
	year = {2018},
	doi = {10.1137/17M1150591}}

@article{HV08,
	author={Hinze, M. and Volkwein, S.},
	title={Error estimates for abstract linear-quadratic optimal control problems using proper orthogonal decomposition},
	journal={Comput. Optim. Appl.},
	volume={39},
	pages={319--345},
	doi={10.1007/s10589-007-9058-4},
	year={2008}}

@book{BMS06,
	title={Dimension Reduction of Large-Scale Systems: Proceedings of a Workshop held in Oberwolfach, Germany, October 19-25, 2003},
	series={Lecture Notes in Computational Science and Engineering},
	year={2005},
	doi = {10.1007/3-540-27909-1},
	publisher={Springer},
	editor = {Benner, P. and Sorensen, D. C. and Mehrmann, V.},
	volume = {45},
	address = {Berlin, Heidelberg}
}

@article{BBH18,
	author = {Baumann, M. and Benner, P. and Heiland, J.},
	title = {Space-time {G}alerkin {POD} with application in optimal control of semi-linear parabolic partial differential equations},
	journal = SIAMSciComp,
	volume = {40},
	pages = {A1611--A1641},
	year = {2018},
	doi = {10.1137/17M1135281}}

@book{QMN16,
	title= {Reduced Basis Methods for Partial Differential Equations: An Introduction},
	author= {Quarteroni, A. and Manzoni, A. and Negri, F.},
	publisher = {Springer},
	series={UNITEXT -- La Matematica per il 3+2},
	address = {Cham},
	doi={10.1007/978-3-319-15431-2},
	year= {2016}}

@book{HRS16,
	title= {Certified Reduced Basis Methods for Parametrized Partial Differential Equations},
	author= {Hesthaven, J.S. and Rozza, G. and Stamm, B.},
	publisher = {Springer},
	series={SpringerBriefs in Mathematics},
	address = {Cham},
	doi = {10.1007/978-3-319-22470-1},
	year= {2016}}

@book{BOCW17,
	title = {Model Reduction and Approximation: Theory and Algorithms},
	publisher = {Society for Industrial and Applied Mathematics},
	year = {2017},
	doi = {10.1137/1.9781611974829},
	address = {Philadelphia, PA},
	editor = {Benner, P. and Ohlberger, M. and Cohen, A. and Willcox, K.},
	series = {Computational Science \& Engineering},
	volume = {15}
}

@article{KarMUV26,
	title={Certified Model Predictive Control for Switched Evolution Equations Using Model Order Reduction},
	author={Michael Kartmann and Mattia Manucci and Benjamin Unger and Stefan Volkwein},
	journal={Commun. Appl. Math. Comput.},
	year={2026},
	doi={10.1007/s42967-025-00562-7}
}

@article{bader2016certified,
	title={Certified reduced basis methods for parametrized distributed elliptic optimal control problems with control constraints},
	author={Bader, Eduard and Kärcher, Mark and Grepl, Martin A and Veroy, Karen},
	doi = {10.1137/16M1059898},
	journal=SIAMSciComp,
	volume={38},
	number={6},
	pages={A3921--A3946},
	year={2016},
	publisher=SIAM
}

@book{bauschke2011convex,
	title={Convex Analysis and Monotone Operator Theory in Hilbert Spaces},
	author={Bauschke, H. H. and Combettes, P. L.},
	isbn={9781441994677},
	lccn={2011926587},
	series={CMS Books in Mathematics},
	year={2011},
	publisher={Springer},
	doi = {10.1007/978-1-4419-9467-7},
	address = {New York}
}

@article{DieG23,
	title = {Reduced order model predictive control for parametrized parabolic partial differential equations},
	author = {Dietze, Saskia and Grepl, Martin A},
	journal = ApplMathComp,
	volume = {453},
	pages = {128044},
	year = {2023},
	publisher = {Elsevier},
	doi = {10.1016/j.amc.2023.128044},
}

@article{azmi2023nonmonotone,
	title = {On the forward--backward method with nonmonotone linesearch for infinite-dimensional nonsmooth nonconvex problems},
	author = {Azmi, Behzad and Bernreuther, Marco},
	journal = {Comput. Optim. Appl.},
	volume = {91},
	number = {3},
	pages = {1263--1308},
	doi = {10.1007/s10589-025-00684-x},
	year = {2025},
}

@article{SchU18,
	author = {Schulze, Philipp and Unger, Benjamin},
	doi = {10.1137/18M1167887},
	journal = SIAMContOpt,
	number = 6,
	pages = {4365--4384},
	title = {Model reduction for linear systems with low-rank switching},
	volume = 56,
	year = 2018
}

@book{troltzsch2024optimal,
	title={Optimal Control of Partial Differential Equations: Theory, Methods and Applications},
	series = {Graduate Studies in Mathematics},
	author={Tr{\"o}ltzsch, Fredi},
	volume={112},
	year={2010},
	publisher={American Mathematical Society},
	doi = {10.1090/gsm/112},
}

@Article{PonGB20,
	author  = {Pontes Duff, I. and Grundel, S. and Benner, P.},
	journal = IEEETransAC,
	title   = {New Gramians for Switched Linear Systems: Reachability, Observability, and Model Reduction},
	year    = {2020},
	number  = {6},
	pages   = {2526--2535},
	volume  = {65},
	doi     = {10.1109/TAC.2019.2934020},
}

@misc{ManU26b,
	author  = {Manucci, M. and Unger, B.},
	title   = {Balancing-based model reduction for switched descriptor systems},
	year    = {2026},
	doi     = {10.48550/arXiv.2404.10511},
	archivePrefix = {arXiv},
	primaryClass = {math.NA}
}

@ARTICLE{MusG91,
	author={Mustafa, D. and Glover, K.},
	journal=IEEETransAC, 
	title={Controller reduction by {$H_\infty$}-balanced truncation},
	year={1991},
	volume={36},
	number={6},
	pages={668--682},
	doi={10.1109/9.86941}}

@Article{PetWL13,
	author   = {Petreczky, M. and Wisniewski, R. and Leth, J.},
	journal  = NonLinAnalHybSys,
	title    = {Balanced truncation for linear switched systems},
	year     = {2013},
	issn     = {1751-570X},
	pages    = {4--20},
	
	volume   = {10},
	doi      = {10.1016/j.nahs.2013.03.007},
}

@inproceedings{PapP14,
	author = {Papadopoulos, A. V. and Prandini, M.},
	booktitle = {Proc. 17th Int. Conf. Hybrid Syst. Comput. Control},
	doi = {10.1145/2562059.2562131},
	pages = {113--122},
	title = {Model Reduction of Switched Affine Systems: {a} Method Based on Balanced Truncation and Randomized Optimization},
	year = {2014},
	address = {Berlin, Germany}
}

@article{PapP16,
	title = {Model reduction of switched affine systems},
	author = {Papadopoulos, A. V. and Prandini, M.}, 
	year = {2016},
	journal = Automatica,
	volume = {70},
	pages = {57--65},
	doi = {10.1016/j.automatica.2016.03.019},
}

@article{WuZ09,
	title = {Weighted $\mathcal{H}_\infty$ model reduction for linear switched systems with time-varying delay},
	author = {Wu, L. and Zheng, W. X.}, 
	year = 2009,
	journal = Automatica,
	volume = 45,
	pages = {186--193},
	doi = {10.1016/j.automatica.2008.06.024},
}

@misc{ManU26,
	title={Model Reduction for Switched Linear Systems via Generalized {Lyapunov} Equations}, 
	author={Mattia Manucci and Benjamin Unger},
	year={2026},
	archivePrefix={arXiv},
	primaryClass={math.NA},
	doi = {10.48550/arXiv.2601.13039}
}

@article{ScaA16,
	author = {Scarciotti, Giordano and Astolfi, Alessandro},
	title = {Model reduction for hybrid systems with state-dependent jumps},
	journal = IFAConline,
	volume = 49,
	number = 18,
	year = 2016,
	pages = {850--855},
	doi = {10.1016/j.ifacol.2016.10.272},
}

@article{BasPWL16,
	author = {Baştuğ, Mert and Petreczky, Mihály and Wisniewski, Rafael and Leth, John},
	journal = IEEETransAC, 
	title = {Model Reduction by Nice Selections for Linear Switched Systems}, 
	year = 2016,
	volume = 61,
	number = 11,
	pages = {3422--3437},
	doi={10.1109/TAC.2016.2518023},
}

@Article{BenGW15,
	author  = {Benner, P. and Gugercin, S. and Willcox, K.},
	journal = {SIAM Rev.},
	title   = {A Survey of Projection-Based Model Reduction Methods for Parametric Dynamical Systems},
	year    = {2015},
	month   = jun,
	pages   = {483--531},
	volume  = {57},
	doi     = {10.1137/130932715},
}

@book{hinrichsen2026mathematical,
	title     = {Mathematical Systems Theory {II}: Control, Observation, Realization, and Feedback},
	author    = {Hinrichsen, Diederich and Pritchard, Anthony J. and Colonius, Fritz and Damm, Tobias and Ilchmann, Achim and Jacob, Birgit and Wirth, Fabian R.},
	year      = {2026},
	publisher = {Springer},
	series = {Texts in Applied Mathematics},
	volume = {85},
	address = {Cham},
	doi = {10.1007/978-3-032-08320-3}
}

@phdthesis{enns1984model,
    author = {Enns, Dale F.},
    title = {Model Reduction for Control System Design},
    school = {Stanford University},
    address = {Stanford, California},
    year = {1984},
}

@inproceedings{breiten2015feedback,
    author = {Breiten, Tobias and Kunisch, Karl},
    title = {Feedback stabilization of the {Schl{\"o}gl} model by {LQG}-balanced truncation},
    booktitle = {2015 European Control Conference (ECC)},
    pages = {1171--1176},
    year = {2015},
    publisher = {IEEE},
    doi = {10.1109/ECC.2015.7330698},
}

@phdthesis{Bader674368,
    author = {Bader, Eduard},
    title = {Reduced Basis Methods Applied to Obstacle Problems and Parametrized Distributed Optimal Control Problems with Control and State Constraints},
    school = {RWTH Aachen University},
    address = {Aachen, Germany},
    year = {2016},
    url = {https://publications.rwth-aachen.de/record/674368},
}

@article{gubisch2014pod,
    author = {Gubisch, Martin and Volkwein, Stefan},
    title = {{POD} a-posteriori error analysis for optimal control problems with mixed control-state constraints},
    journal = {Comput. Optim. Appl.},
    volume = {58},
    number = {3},
    pages = {619--644},
    year = {2014},
    doi = {10.1007/s10589-014-9636-1},
}

@STRING{ApplMathComp = {Appl. Math. Comput.}}

@STRING{Automatica = {Automatica}}

@STRING{IEEE = {{IEEE}}}

@STRING{IEEETransAC = {{IEEE} Trans. Automat. Control}}

@STRING{IFAConline = {IFAC-PapersOnLine}}

@STRING{NonLinAnalHybSys = {Nonlinear Anal. Hybrid Syst.}}

@STRING{SCL = {Syst. Control Lett.}}

@STRING{SIAM = {Society for Industrial and Applied Mathematics}}

@STRING{SIAMContOpt = {{SIAM} J. Control Optim.}}

@STRING{SIAMSciComp = {{SIAM} J. Sci. Comput.}}

@STRING{Springer = {Springer-Verlag}}

\end{document}